\documentclass{article}
\usepackage{amsmath,amssymb,amsthm}
\def\Cmath{\mathbb{C}}
\def\Nmath{\mathbb{N}}
\def\Rmath{\mathbb{R}}
\def\Smath{\mathbb{S}}
\def\Zmath{\mathbb{Z}}
\newcommand{\Aut}[1]{Aut\,{#1}}
\newcommand{\bv}[1]{\mathbf{#1}}
\newtheorem{lemma}{Lemma}
\newtheorem{remark}{Remark}
\newtheorem{example}{Example}
\newtheorem{theorem}{Theorem}
\newtheorem{corollary}{Corollary}
\title{Left invariant complex structures on $U(2)$ and $SU(2) \times SU(2)$ revisited.
\thanks{\textit{Math. Subj. Class. [2000]}  : 17B30 }
}
\author{L. Magnin \\
Institut Math\'{e}matique de Bourgogne
\thanks{UMR CNRS 5584, Universit\'{e} de Bourgogne, BP 47870, 21078 Dijon Cedex, France.}
\\
\texttt{magnin@u-bourgogne.fr}
}
\date{\today}
\begin{document}
\maketitle
\begin{abstract}
We compute the torsion-free linear maps
$J \, : \, \mathfrak{su}(2) \rightarrow \mathfrak{su}(2),$
deduce a new determination of the
integrable complex structures
and their equivalence classes under the action of the automorphism group
for $\frak{u}(2)$ and
$\mathfrak{su}(2) \oplus \mathfrak{su}(2),$
and prove
that in both cases
the set of complex structures is a differentiable manifold.
$\mathfrak{u}(2) \oplus \mathfrak{u}(2),$
$\mathfrak{su}(2)^N$ and $\mathfrak{u}(2)^N$ are also considered.
Extension of complex structures from
$\frak{u}(2)$
to $\mathfrak{su}(2) \oplus \mathfrak{su}(2)$
are  studied,
local holomorphic charts given,
and attention is paid to what representations of
$\frak{u}(2)$ we can get from a substitute to the regular representation
on a space of holomorphic functions for the complex structure.
\end{abstract}
\section{Introduction.}
The left invariant integrable complex structures on the group
$U(2)$ of unitary $2\times 2$ matrices,
 i.e. integrable complex structures on its Lie algebra $\frak{u}(2),$
have been computed for the first time, up to equivalence, in \cite{sasaki}
in the algebraic approach, that is
by determining the complex Lie subalgebras
$\frak{m}$
of the complexification
of
$\mathfrak{g}_{\Cmath}$ of
$\mathfrak{u}(2)$
such that
$\mathfrak{g}_{\Cmath} = \mathfrak{m} \oplus \bar{\mathfrak{m}},$
bar denoting conjugation.
More recently, and more generally,
all left invariant maximal rank $CR$-structures
on any finite dimensional compact Lie group have been classified up to equivalence
in \cite{charbonnel}.
Independently, in the case of  $SU(2)\times SU(2),$ the integrable complex structures on
$\mathfrak{su}(2) \oplus \mathfrak{su}(2)$ have been computed
in \cite{daurtseva} by direct approach and computations.

In the present paper,
we  first compute  the torsion-free linear maps
$J \, : \, \mathfrak{su}(2) \rightarrow \mathfrak{su}(2).$
They appear to be maximal rank $CR$-structures, of the $CR0$-type in
the classification of \cite{charbonnel}.
Then we show how to deduce,
 with the computer assisted methods of \cite{artmagnin1},
a new determination of the
integrable complex structures
and their equivalence classes under the action of the automorphism group
for  the specific cases of $\frak{u}(2)$
and
$\mathfrak{su}(2) \oplus \mathfrak{su}(2),$
without resorting to the general results of \cite{charbonnel}.
Our method  consists in growing dimensions starting with
torsion-free linear maps of
$\mathfrak{su}(2).$

$\mathfrak{u}(2) \oplus \mathfrak{u}(2),$
and extensions for $\mathfrak{su}(2)^N$ and $\mathfrak{u}(2)^N$ are done too.
In thoses cases,
the set of complex structures is a differentiable manifold,
though  we write down  explicit proofs only
in the cases of
$\mathfrak{u}(2)$ and
$\mathfrak{su}(2) \oplus \mathfrak{su}(2).$
We also
examine the extension of complex structures from
$\frak{u}(2)$
to $\mathfrak{su}(2) \oplus \mathfrak{su}(2),$
compute  local complex charts
for the complex manifolds associated to the complex structures,
and determine what representations of
$\frak{u}(2)$ we can get from a substitute to the regular representation
on a space of holomorphic functions for the complex structure.

\section{Preliminaries.}
Let $G_0$ be a connected finite dimensional real Lie group, with
Lie algebra $\mathfrak{g}.$
An almost complex structure on
$\mathfrak{g}$ is a linear map $J \, : \, \mathfrak{g} \rightarrow  \mathfrak{g}$
such that $J^2=-1.$
The almost complex structure $J$ is said to be \textit{integrable} if it satisfies the condition
\begin{equation}
\label{Jinv}
 [{J} X, {J}Y]-[X,Y]-{J}[{J}X,Y]-J[X,{J}Y] =0 \quad \forall X,Y \in \mathfrak{g}.
 \end{equation}
From the Newlander-Nirenberg theorem \cite{NN}, condition (\ref{Jinv})
means that $G_0$
can be given the structure of a complex manifold with the same underlying
real structure and such that the canonical complex structure on $G_0$ is the left invariant almost
complex structure
$\hat{J}$ associated to $J.$
(For more details, see \cite{artmagnin1}).
By a complex structure  on $\mathfrak{g}$, we will mean an
\textit{integrable} almost complex structure on  $\mathfrak{g},$ that is one satisfying
(\ref{Jinv}).

\par
Let $J$ a complex structure on $\mathfrak{g}$ and denote by $G=(G_0,J)$
the group $G_0$ endowed with the structure
of complex manifold  defined by $\hat{J}.$
The complexification $\mathfrak{g}_{\Cmath}$ of $\mathfrak{g}$ splits
as $\mathfrak{g}_{\Cmath}= \mathfrak{g}^{(1,0)} \oplus \mathfrak{g}^{(0,1)}$
where $\mathfrak{g}^{(1,0)}= \{\tilde{X}=X-iJX; \, X \in \mathfrak{g}\},$
$\mathfrak{g}^{(0,1)}=\{
\tilde{X}^-=
X+iJX; \, X \in \mathfrak{g}\}.$
We will denote  $\mathfrak{g}^{(1,0)}$ by $\mathfrak{m}.$ The integrability of $J$
amounts to
$\mathfrak{m}$ being a complex subalgebra of
$\mathfrak{g}_{\Cmath}.$ In that way  the set of complex structures on $\mathfrak{g}$ can be identified with the set of all complex subalgebras $\mathfrak{m}$ of
$\mathfrak{g}_{\Cmath}$ such that
$\mathfrak{g}_{\Cmath} = \mathfrak{m} \oplus \bar{\mathfrak{m}},$  bar denoting conjugation
in $\mathfrak{g}_{\Cmath}.$ In particular, $J$ is said to be abelian if $\mathfrak{m}$ is.
That is the algebraic approach.
Our approach  is more trivial. We fix a basis of $\mathfrak{g},$
write down
the torsion equations
$ij|k$ ($1 \leqslant i,j,k \leqslant n$ ) obtained by projecting on $x_k$ the equation
 $[{J} x_i, {J}x_j]-[x_i,x_j]-{J}[{J}x_i,x_j]-J[x_i,{J}x_j] =0,$
 where
$(x_j)_{ 1 \leqslant j \leqslant n}$ is the basis of $\mathfrak{g}$ we use,
and solve them in steps
by specific programs with the computer algebra system \textit{Reduce} by A. Hearn.
These programs, are downloadable in
\cite{companionarchive}.
From now on, we will use the same notation $J$ for $J$ and $\hat{J}$ as well.
For any $x \in G_0,$
the complexification ${T_x(G_0)}_{\Cmath}$ of the tangent space
also splits as the direct sum of the holomorphic vectors
${T_x(G_0)}^{(1,0)}= \{
\tilde{X}=
X-iJX; \, X \in T_x(G_0)\}$ and the antiholomorphic vectors
${T_x(G_0)}^{(0,1)}= \{
\tilde{X}^-=
X+iJX; \, X \in T_x(G_0)\}.$
For any open subset $V \subset G_0,$ the space  $H_{\Cmath}(V)$  of complex valued holomorphic functions on $V$
is comprised of all complex smooth functions $f$ on $V$
 which are annihilated by any antiholomorphic vector field. This is equivalent
to $f$ being annihilated by
all
\begin{equation}
\label{Xj-tilde}
\tilde{X}_j^{-} = X_j +i J X_j   \quad 1 \leqslant j \leqslant n
\end{equation}
with $(X_j)_{1 \leqslant j \leqslant n}$ the left invariant vector fields associated
to the basis $(x_j)_{1 \leqslant j \leqslant n}$ of
$\mathfrak{g}.$
Hence :
\begin{equation}
H_{\Cmath}(V) =\{f \in C^{\infty}(V) \; ; \; \tilde{X}_j^{-} \, f = 0
\; \; \forall j\;\;  1 \leqslant j \leqslant n\}.
\end{equation}
\par
Finally, the automorphism group
$\text{Aut } \mathfrak{g}$ of $\mathfrak{g}$ acts on the set
$\mathfrak{X}_{\mathfrak{g}}$ of all complex structures on  $\mathfrak{g}$ by
$J \mapsto \Phi \circ J \circ \Phi^{-1} \quad \forall \Phi \in
\text{Aut } \mathfrak{g}.$
Two complex structures $J,J^{\prime}$
on   $\mathfrak{g}$  are said to be \textit{equivalent} if they are on the same
$\text{Aut } \mathfrak{g}$ orbit.
For simply connected $G_0,$
this amounts to requiring the existence of an $f \in \text{Aut } G_0$ such that
$f : (G_0,J) \rightarrow  (G_0,J^{\prime})$ is biholomorphic.

\section{$U(2)$.}
Consider the Lie algebra $\mathfrak{su}(2)$ along with its basis $\{J_1,J_2,J_3\}$
defined by
$J_1 =\frac{i}{2} \left(\begin{smallmatrix} 0&1 \\ 1&0  \end{smallmatrix}\right)$,
$J_2 =\frac{1}{2} \left(\begin{smallmatrix} 0&-1 \\ 1&0 \end{smallmatrix}\right)$,
$J_3 =\frac{i}{2} \left(\begin{smallmatrix} 1&0 \\ 0&-1 \end{smallmatrix}\right).$ One has
\begin{equation}
\label{relationssu2}
[J_1,J_2]=J_3, \; [J_2,J_3]=J_1, \; [J_3,J_1]=J_2
\end{equation}
and the corresponding one-parameter subgroups of $SU(2)$ are
$e^{tJ_1} =
\left(\begin{smallmatrix}
\cos{\frac{t}{2}}
&i\sin{\frac{t}{2}}\\
i\sin{\frac{t}{2}}
&\cos{\frac{t}{2}}\\
\end{smallmatrix}\right)$,
$e^{tJ_2} =
\left( \begin{smallmatrix}
\cos{\frac{t}{2}}
&-\sin{\frac{t}{2}}\\
\sin{\frac{t}{2}}
&\cos{\frac{t}{2}}\\
\end{smallmatrix} \right)$,
$e^{tJ_3} =
\left( \begin{smallmatrix}
e^{i\frac{t}{2}}
&0\\ 0&
e^{-i\frac{t}{2}}
\end{smallmatrix} \right)$.
By means of the basis
$\{J_1,J_2,J_3\}$,
$\mathfrak{su}(2)$ can be identified to the euclidean vector space $\Rmath^3$
the bracket being then identified to the vector
product $\wedge$. Then
$\Aut{\mathfrak{su}(2)}$ is comprised of the matrices
$A = Mat ( \bf{a},\bf{b}, \bf{a}\wedge \bf{b})$ with $\bf{a},\bf{b}$  any two
orthogonal normed vectors in
$\Rmath^3$, i.e.
$\Aut{\mathfrak{su}(2)}\cong SO(3)$.
Now,
$\mathfrak{u}(2) = \mathfrak{su}(2) \oplus \mathfrak{c}$
where $\mathfrak{c}= \Rmath \, J_4$ is the center of
$\mathfrak{u}(2),$
$J_4 = \frac{i}{2}\left(\begin{smallmatrix} 1&0 \\ 0&1 \end{smallmatrix}\right).$
We use the basis $(J_1,J_2,J_3,J_4)$ for
$\mathfrak{u}(2).$
$\Rmath^*$ stands for $\Rmath \setminus \{ 0\}.$
\begin{lemma}
$\Aut{\mathfrak{u}(2)}\cong SO(3) \times \Rmath^*.$
\end{lemma}
\begin{proof}
As the center is invariant, any
$\Phi \in \Aut{\mathfrak{u}(2)}$ is of the form
$$\Phi = \begin{pmatrix} &A&& \begin{matrix} 0\\0\\0 \end{matrix}\\ b^4_1&b^4_2&b^4_3&b^4_4\end{pmatrix}$$
with $A \in  \Aut{\mathfrak{su}(2)}\cong SO(3)$ and $b^4_4 \in \Rmath^*.$
Necessarily,
$b^4_1 =b^4_2=b^4_3=0,$
since $\Phi(J_k) \in [\mathfrak{u}(2),\mathfrak{u}(2)]=\mathfrak{su}(2)$ ($1\leqslant k \leqslant 3$).
\end{proof}

\begin{lemma}
\label{su2notorsion}
Let  $J \, : \, \mathfrak{su}(2)  \rightarrow  \mathfrak{su}(2)$ linear.
 $J$ has zero torsion, i.e. satisfies
(\ref{Jinv}), if and only if
there exists  $R \in  SO(3)$
such that
\begin{equation}
\label{Jsu2notorsion}
R^{-1}J R =
\begin{pmatrix}
0&1&0\\
-1&0&0\\
0&0&\xi^3_3
\end{pmatrix}.
\end{equation}
\end{lemma}
\begin{proof}
Let $J = (\xi^i_j)_{1 \leqslant i,j \leqslant 3}$
in the basis ($J_1,J_2,J_3$).
The 9 torsion equations are
\begin{eqnarray*}
12|1& \quad \xi^1_3( \xi^2_2 +\xi^1_1) + \xi^3_1( \xi^2_2 -\xi^1_1) - \xi^3_2(\xi^2_1 +\xi^1_2) =0,\\
12|2& \quad \xi^2_3( \xi^2_2 +\xi^1_1) - \xi^3_2( \xi^2_2 -\xi^1_1) - \xi^3_1(\xi^2_1 +\xi^1_2) =0,\\
12|3& \quad \xi^1_2 \xi^2_1 -\xi^2_2\xi^1_1 -(\xi^3_1)^2 -(\xi^3_2)^2 +\xi^3_3(\xi^2_2+\xi^1_1) +1 =0,\\
13|1& \quad \xi^1_1( \xi^2_1 -\xi^1_2) + \xi^2_3(\xi^1_3 +\xi^3_1) - \xi^3_3(\xi^2_1+\xi^1_2) =0,\\
13|2& \quad \xi^1_3 \xi^3_1 + \xi^2_2\xi^1_1 -(\xi^2_1)^2 -(\xi^2_3)^2 +\xi^3_3(\xi^2_2-\xi^1_1) +1 =0,\\
13|3& \quad -\xi^1_1( \xi^2_3 +\xi^3_2) + \xi^2_1( \xi^1_3 +\xi^3_1) + \xi^3_3(\xi^2_3 -\xi^3_2) =0,\\
23|1& \quad \xi^3_2 \xi^2_3 +\xi^2_2\xi^1_1 -(\xi^1_3)^2 -(\xi^1_2)^2 -\xi^3_3(\xi^2_2-\xi^1_1) +1 =0,\\
23|2& \quad \xi^2_2( \xi^2_1 -\xi^1_2) - \xi^1_3( \xi^2_3 +\xi^3_2) + \xi^3_3(\xi^2_1 +\xi^1_2) =0,\\
23|3& \quad \xi^2_2( \xi^3_1 + \xi^1_3) - \xi^1_2( \xi^2_3 +\xi^3_2) + \xi^3_3(\xi^3_1 -\xi^1_3) =0.
\end{eqnarray*}

Again, we identify  $\mathfrak{su}(2)$ to $\Rmath^3$ with the vector product
by means of  the basis ($J_1,J_2,J_3$).
$J$ has at least one  real eigenvalue $\lambda.$
Let  $\bf{f_3} \in \Rmath^3$
some normed eigenvector associated to $\lambda$. Then there exist normed vectors
$\bf{f_1},\bf{f_2} \in \Rmath^3$ such that
$(\bf{f_1},\bf{f_2},\bf{f_3})$ is a direct orthonormal basis of  $\Rmath^3.$
Hence there exists $R \in SO(3)$
such that
$$R^{-1} A R = \begin{pmatrix} * & * & 0\\ * &* & 0\\ *&*&\lambda \end{pmatrix}.$$
Hence we may suppose $\xi^1_3 = \xi^2_3$ in $J.$
Now, the torsion equations $12|1$ and  $12|2$ read respectively
$\xi^3_1 ( \xi^2_2  - \xi^1_1) = \xi^3_2 ( \xi^2_1  + \xi^1_2),$
$\xi^3_1 ( \xi^2_1  + \xi^1_2) = \xi^3_2 ( \xi^1_1  - \xi^2_2)
$
and imply the 2 equations
$(\xi^3_1)^2 ( \xi^2_2  - \xi^1_1) = -(\xi^3_2)^2 ( \xi^2_2  - \xi^1_1),$
$(\xi^3_2)^2 ( \xi^2_1  + \xi^1_2) = -(\xi^3_1)^2 ( \xi^2_1  + \xi^1_2).$
Hence each one of the conditions $\xi^2_2 \neq \xi^1_1 $ or
$\xi^2_1 \neq -\xi^1_2 $ implies
$\xi^3_1 = \xi^3_2 =0.$
We now have 2 cases. Case 1 : $\xi^3_1 = \xi^3_2 =0.$  Case 2 :
$\xi^3_1, \xi^3_2$ not both zero.
In Case 2, one necessarily has
$\xi^2_2 = \xi^1_1 $ and
$\xi^2_1 = -\xi^1_2 .$ Then equations $23|1$ and $23|2$ read
$-(\xi^1_2)^2 + (\xi^1_1)^2 +1 =0,$ $\xi^1_2 \xi^1_1 =0$ and give $\xi^1_1=0$,
$\xi^1_2 =\pm 1.$ Now equation $12|3$ reads $(\xi^3_2)^2 + (\xi^3_1)^2 =0.$
Hence Case 2 doesn't occur, i.e. one may suppose $\xi^3_1=\xi^3_2=0.$
Then equations  $13|1$ and $23|2$ read resp.
$\xi^3_3 (\xi^2_1 +\xi^1_2) = \xi^1_1 (\xi^2_1 -\xi^1_2),$
$\xi^3_3 (\xi^2_1 +\xi^1_2) = -\xi^2_2 (\xi^2_1 -\xi^1_2),$
hence if $\xi^2_1 \neq \xi^1_2,$ necessarily $\xi^2_2= -\xi^1_1.$
Now  $\xi^2_1 = \xi^1_2$ is impossible since it would imply
either $\xi^3_3=0$ or $\xi^1_2=0.$ In fact, first, if $\xi^1_2=0$,
 equations  $12|3,$ $13|2$ and  $23|1$ read resp.
$\xi^3_3 (\xi^2_2 +\xi^1_1) - \xi^2_2 \xi^1_1 +1 =0,$
$\xi^3_3 (-\xi^2_2 +\xi^1_1) - \xi^2_2 \xi^1_1 -1 =0,$
$\xi^3_3 (-\xi^2_2 +\xi^1_1) + \xi^2_2 \xi^1_1 +1 =0,$
so that $12|3 + 13|2$ gives $\xi^1_1(\xi^3_3 -\xi^2_2)=0$ and
$12|3 + 23|1$ gives $\xi^1_1\xi^3_3 =-1,$ hence $\xi^1_1 \neq 0$ and $\xi^3_3= \xi^2_2,$
which is impossible since then $12|3$ reads $(\xi^2_2)^2+1=0.$ Second, if
$\xi^3_3=0,$
$12|3$, $13|2$
read resp.
$-\xi^2_2 \xi^1_1 +(\xi^1_2)^2 +1 =0,$
$-\xi^2_2 \xi^1_1 +(\xi^1_2)^2 -1 =0,$
which is contradictory.
Hence we get as asserted $\xi^2_1 \neq \xi^1_2$ and $\xi^2_2= -\xi^1_1.$
Now we prove that $\xi^1_1=0$ and $\xi^1_2=\pm 1.$
Since $\xi^2_2 =- \xi^1_1,$
equations $12|3,$ $13|1,$ $13|2,$ $ 23|1$
read resp.
\begin{eqnarray*}
12|3& \quad \xi^1_2 \xi^2_1 + (\xi^1_1)^2 +1 =0,\\
13|1& \quad \xi^3_3( \xi^2_1 +\xi^1_2) - \xi^1_1(\xi^2_1 -\xi^1_2) =0,\\
13|2& \quad 2\xi^3_3 \xi^1_1 +(\xi^2_1)^2 + (\xi^1_1)^2-1 =0,\\
23|1& \quad 2\xi^3_3 \xi^1_1 -(\xi^1_2)^2 - (\xi^1_1)^2+1 =0.\\
\end{eqnarray*}
From $12|3,$ $\xi^1_2 \neq 0$ and
$ \xi^2_1=-\frac{1+(\xi^1_1)^2}{\xi^1_2}.$
Then $13|1,$
 $13|2,$
read resp. $Q=0,$ $R=0$ with
$Q= \xi^1_1( (\xi^1_1)^2 +(\xi^1_2)^2 +1) -\xi^3_3( (\xi^1_1)^2 -(\xi^1_2)^2 +1),$
$R= (\xi^1_2)^2 (2\xi^3_3 \xi^1_1 +(\xi^1_1)^2 -1) +( (\xi^1_1)^2 +1)^2.$
Denote from
$23|1$ $S= 2\xi^3_3 \xi^1_1 -(\xi^1_2)^2 - (\xi^1_1)^2+1.$
Suppose $\xi^1_1 \neq 0.$ Then
$N=\frac{R-S}{\xi^1_1}=
2\xi^3_3 ((\xi^1_2)^2 -1)  +\xi^1_1 ( (\xi^1_1)^2 +(\xi^1_2)^2 +3)=0$
would give, for $\xi^1_2 \neq  \pm 1,$
$\xi^3_3= - \frac{
\xi^1_1 ( (\xi^1_1)^2 +(\xi^1_2)^2 +3)
}{2((\xi^1_2)^2-1)}$
and then
$R=
- \frac{(
(\xi^1_2)^2 -2\xi^1_2 +(\xi^1_1)^2 +1)
((\xi^1_2)^2 +2\xi^1_2 +(\xi^1_1)^2 +1)
}{(\xi^1_2)^2-1}$
which is impossible since the polynomial $X^2 \pm 2X + (\xi^1_1)^2 +1$
has no real root.
Hence $\xi^1_2= \pm 1.$
Now $S=0$ gives $\xi^3_3=\frac{\xi^1_1}{2}$  and then $R=
(\xi^1_1)^2 ((\xi^1_1)^2 +4) \neq 0.$
Hence $\xi^1_1=0.$
Finally that implies as asserted $\xi^1_2=\pm 1,$ since $R=-(\xi^1_2)^2+1.$
We conclude that $\xi^1_1=\xi^2_2=0, \xi^2_1=-\xi^1_2, \xi^1_2=\pm 1.$
Changing if necessary $\Phi$ to $\Phi\,  \text{diag}(\left(\begin{smallmatrix} 0&1\\1&0\end{smallmatrix} \right),-1)$), one may suppose $\xi^1_2=1.$
\end{proof}
\begin{remark}
\rm
Recall that a rank $r$ $CR$-structure on a real Lie algebra $\frak{g}$ is  a
$r$-dimensional subalgebra $\frak{m}$ of the complexification $\frak{g}_\Cmath$ of $\frak{g}$
such that $\frak{m} \cap \bar{\frak{m}} =\{0\}.$
Then $\frak{m}=\{X-iJ_\frak{p}X; X\in\frak{p}\}$
where $\frak{p}$ (the real part of $\frak{m}$) is a vector subspace of $\frak{g}$ and
$J_\frak{p} \, : \, \mathfrak{p}  \rightarrow  \mathfrak{p}$ is a   zero torsion linear map
such that $J_\frak{p}^2= -I\negthinspace d_{\frak{p}}.$
Alternatively, a $CR$-structure can be defined by such data $(\frak{p}, J_\frak{p}).$
For even-dimensional $\frak{g},$
$CR$-structures of maximal rank $r=\frac{1}{2}\, \dim{\frak{g}}$ are just complex structures on $\frak{g}.$
$CR$-structures of maximal rank on a real compact Lie algebra have been classified in \cite{charbonnel}.
For odd-dimensional $\frak{g},$
they fall  essentially into 2 classes: $CR0$ and (strict)  $CRI$.
For  even-dimensional $\frak{g}$ they are all $CR0.$
From lemma \ref{su2notorsion},
any linear map  $J \, : \, \mathfrak{su}(2)  \rightarrow  \mathfrak{su}(2)$
which  has zero torsion is such that $\ker{(J^2+I\negthinspace d)} \neq \{0\},$
 and hence defines a maximal rank $CR$-structure on
$\frak{su}(2).$
It is of type $CR0.$ Let us elaborate on that point.
$\frak{a}_0= \Cmath J_3$ is a Cartan subalgebra of
$\frak{su}(2).$
The complexification $\frak{sl}(2)$ of $\frak{su}(2)$ decomposes as
$\frak{sl}(2)= \Cmath H_- \oplus \frak{h} \oplus \Cmath H_+$
 with
 $H_{\pm} =iJ_1 \mp J_2, \; H_3=iJ_3 ,\; \frak{h}= \Cmath H_3.$ Any maximal rank  $CR$-structure
 of  $CR0$-type (resp. (strict) $CRI$-type)
 is equivalent to $\frak{m} = \Cmath H_+$ (resp. $\frak{m}= \Cmath(a J_3 +H_+)$, $a \in \Rmath^*$), and has
 real part $\frak{p}= \Rmath J_1 \oplus \Rmath J_2$
 (resp.
 $\frak{p}= \Rmath J_1 \oplus \Rmath J^{\prime}_2, \; J^{\prime}_2= J_2-aJ_3$).
 The corresponding endomorphism $J_\frak{p}$ of $\frak{p}$
 has matrix
 $\left(\begin{smallmatrix}0&1\\-1&0 \end{smallmatrix} \right)$
 in the basis $(J_1,J_2)$ (resp. $(J_1,J^{\prime}_2)$)
 and has zero torsion on $\frak{p}.$
 Any extension of
 $J_\frak{p}$ to
$\frak{su}(2)$ has matrix
 $\left(\begin{smallmatrix}0&1&\xi^1_3\\-1&0&\xi^2_3\\0&0&\xi^3_3 \end{smallmatrix} \right)$
 in the basis  $(J_1,J_2,J_3)$ (resp. $(J_1,J^{\prime}_2,J_3)$).
 In the $CR0$ case, it   has zero torsion on the whole of  $\frak{su}(2)$
 if and only if $\xi^1_3=\xi^2_3=0,$ i.e. is of the form
(\ref{Jsu2notorsion}).
 In the $CRI$ case, it never has
 zero torsion on the whole of  $\frak{su}(2).$
\end{remark}

\begin{lemma}
\label{product}
Let $\frak{g} = \bigoplus_{j=1}^N \, \frak{g}^{(j)},$
where $\frak{g}^{(j)}$ are real Lie algebras
with bases $\mathcal{B}_j =(X^{(j)}_{k})_{1\leqslant k\leqslant n_j},$
and   let $\pi^{(j)}  \, : \, \frak{g}  \rightarrow  \frak{g}^{(j)}$ be the projections.
Let  $J \, : \, \frak{g}  \rightarrow  \frak{g}$ be a linear map,
 $\pi^i_j = \pi^{(i)}\circ J \circ \pi^{(j)},$
 $\tilde{\pi}^i_j = \pi^{(i)}\circ J \circ \pi^{(j)}|_{\frak{g}^{(j)}}.$
If $J$ has zero torsion, then the 2 following conditions are satisfied:
\\
(i) $\tilde{\pi}^i_i$ has zero torsion for any $i;$
\\ (ii)
$[\pi^i_j X,\pi^i_j Y]=
\pi^i_j [JX,Y] + \pi^i_j [X,JY]\; \forall X,Y \in \frak{g}^{(j)}$
for any $i,j$ such that $i\neq j.$
\end{lemma}
\begin{proof}
 For any $i,j$ let $X,Y \in \frak{g}.$
Applying $\pi^{(i)}$ to the torsion equation
(\ref{Jinv})
we get
\begin{equation}
\label{eqlemme}
[\pi^{(i)}JX,\pi^{(i)}JY]
-[\pi^{(i)}X,\pi^{(i)}Y]
-\pi^{(i)}J[JX,Y]
-\pi^{(i)}J[X,JY]=0.
\end{equation}
Suppose first $i=j$ and
 $X,Y \in \frak{g}^{(i)}.$  Then
$[JX,Y]=
[\pi^{(i)}JX,Y]
=\pi^{(i)}[\pi^{(i)}J\pi^{(i)}X,Y],$
and $[X,JY]=
[X,\pi^{(i)}JY]
=\pi^{(i)}[X,\pi^{(i)}J\pi^{(i)}Y],$
and moreover
$[\pi^{(i)}X,\pi^{(i)}Y]
=[X,Y]$,
hence (\ref{eqlemme}) gives
$$[\pi^{(i)}JX,\pi^{(i)}JY]
-[X,Y]
-\pi^{(i)}J\pi^{(i)}[\pi^{(i)}J\pi^{(i)}X,Y]
-\pi^{(i)}J\pi^{(i)}[X,\pi^{(i)}J\pi^{(i)}Y],$$
i.e.
$$[\tilde{\pi}^i_i X,\tilde{\pi}^i_iY]
-[X,Y]
-\tilde{\pi}^i_i[\tilde{\pi}^i_iX,Y]
-\tilde{\pi}^i_i[X,\tilde{\pi}^i_iY]$$
that is $\tilde{\pi}^i_i$ has no torsion.
Suppose now $i\neq j$
and  $X,Y \in \frak{g}^{(j)}.$
Then $[\pi^{(i)}X,\pi^{(i)}Y]=0$
and (\ref{eqlemme}) gives
$$
[\pi^{(i)}JX,\pi^{(i)}JY]
-\pi^{(i)}J\pi^{(j)}[JX,Y]
-\pi^{(i)}J\pi^{(j)}[X,JY]=0,
$$
i.e.
$$[{\pi}^i_j X,{\pi}^i_jY]
-{\pi}^i_j[JX,Y]
-{\pi}^i_j[X,JY]=0.$$
\end{proof}

\begin{theorem}
\label{theoremu2notorsion}
(i) Let  $J \, : \, {\mathfrak{u}(2) \rightarrow
 \mathfrak{u}(2)}$ linear.
 $J$ has zero torsion, i.e. satisfies
(\ref{Jinv}), if and only if
there exists
$\Phi \in  SO(3)\times \Rmath_+^*$
such that
\begin{equation}
\label{Ju2notorsion}
\Phi^{-1}J\Phi =
\begin{pmatrix}
0&1&0&0\\
-1&0&0&0\\
0&0&\xi^3_3&\xi^3_4\\
0&0&\xi^4_3&\xi^4_4
\end{pmatrix},
\quad
\begin{pmatrix}\xi^3_3&\xi^3_4\\
\xi^4_3&\xi^4_4
\end{pmatrix} \in \textit{gl}(2,\Rmath).
\end{equation}
\\ (ii) Any $J \in {\frak{X}}_{\mathfrak{u}(2)}$
is equivalent to a unique
\begin{equation}
\label{Ju2m}
J(\xi) =
\begin{pmatrix}
0&1&0&0\\
-1&0&0&0\\
0&0&\xi&1\\
0&0&-(1+{\xi}^2)&-\xi
\end{pmatrix}
\end{equation}
with $\xi\in \Rmath.$
$J(\xi)$ and $J(\xi^{\prime})$  ($\xi,\xi^{\prime} \in \Rmath$)
are equivalent
if and only if
$\xi=\xi^{\prime}.$
\end{theorem}
\begin{proof}
(i)
From lemma \ref{product}
$$J=\begin{pmatrix}
&&&\xi^1_4\\
&J_1&&\xi^2_4\\
&&&\xi^3_4\\
\xi^4_1&\xi^4_2&\xi^4_3&\xi^4_4
\end{pmatrix}
$$
for some
  $J_1 \, : \, \mathfrak{su}(2)  \rightarrow  \mathfrak{su}(2)$
  with zero torsion.
  From lemma \ref{su2notorsion},
there exists  $R \in  SO(3)$
such that
$R^{-1}J_1 R =
\left(\begin{smallmatrix}
0&1&0\\
-1&0&0\\
0&0&\xi^3_3
\end{smallmatrix}\right)$ whence
$$\Phi^{-1}J\Phi =
\begin{pmatrix}
0&1&0&\xi^1_4\\
-1&0&0&\xi^2_4\\
0&0&\xi^3_3&\xi^3_4\\
\xi^4_1&\xi^4_2&\xi^4_3&\xi^4_4
\end{pmatrix}$$
$\Phi= \text{diag}(R,1).$
Hence we may suppose
$J_1  =
\left(\begin{smallmatrix}
0&1&0\\
-1&0&0\\
0&0&\xi^3_3
\end{smallmatrix}\right).$
Now      the  torsion equations
$13|4,23|4$
$14|3,24|3$ give the 2 Cramer systems
$\xi^4_2\xi^3_3 +\xi^4_1=0,$
$-\xi^4_2 + \xi^3_3 \xi^4_1=0;$
$\xi^2_4\xi^3_3 -\xi^1_4=0,$
$\xi^2_4 + \xi^3_3 \xi^1_4=0.$
Hence $\xi^4_1=\xi^4_2=\xi^1_4=\xi^2_4=0.$
Then all torsion equations vanish, and  (i) is proved
(\cite{companionarchive}, $\texttt{torsionu2.red}$).
\\(ii)
From (i), we may suppose
$$J=\begin{pmatrix}
0&1&0&0\\
-1&0&0&0\\
0&0&\xi^3_3&\xi^3_4\\
0&0&\xi^4_3&\xi^4_4
\end{pmatrix}.$$
Now
$J \in {\frak{X}}_{\mathfrak{u}(2)}$
if and only if
$\left(\begin{smallmatrix}\xi^3_3&\xi^3_4\\
\xi^4_3&\xi^4_4
\end{smallmatrix} \right)^2 = -I,$
i.e.
$$J=\begin{pmatrix}
0&1&0&0\\
-1&0&0&0\\
0&0&\xi^3_3&\xi^3_4\\
0&0&-\frac{1+(\xi^3_3)^2}{\xi^3_4}&-\xi^3_3
\end{pmatrix},\quad \xi^3_4 \neq 0.$$
Observe now that for any
$\Phi= \text{diag}(A,b) \in \Aut{\mathfrak{u}(2)} \quad (A \in SO(3), b \neq 0)$
\begin{equation}
\label{equivu2}
\Phi J \Phi^{-1} =
\begin{pmatrix} &AJ_1A^{-1}&& b^{-1}A
\begin{pmatrix} 0\\ 0 \\ \xi^3_4 \end{pmatrix}
\\
&b
\begin{pmatrix} 0& 0 &
-\frac{1+(\xi^3_3)^2}{\xi^3_4}
\end{pmatrix}A^{-1}
&&-\xi^3_3\end{pmatrix}.
\end{equation}
Taking $A=I,$ $b= \xi^3_4,$
we get
$$
\Phi J \Phi^{-1} =
\begin{pmatrix}
0 & 1 &0&0\\
-1 & 0 &0 &0\\
0 & 0 &\xi^3_3 & 1 \\
0 & 0 &-(1+(\xi^3_3)^2)  &-\xi^3_3
\end{pmatrix}
.
$$
Hence $J$ is equivalent  to
$J(\xi)$ in (\ref{Ju2m}) with $\xi =\xi^3_3.$
The last assertion of the theorem results from (\ref{equivu2}).
\end{proof}

\begin{remark}
\rm
In \cite{sasaki}, the equivalence classes of left invariant integrable complex structures
on
${\mathfrak{u}(2)}$
are shown to be parametrized by the complex subalgebras
$\frak{m}_d$ with basis
$\{J_1+iJ_2,2iJ_3+dJ_4\}$ with $d=-\frac{1+i\xi}{1+\xi^2}, \, \xi \in \Rmath.$
The complex structure defined by
$\frak{m}_d$ has matrix
\begin{equation*}
\begin{pmatrix}
0 &1& 0&0\\
-1 &0& 0&0\\
0 &0& \xi&2(1+\xi^2)\\
0 &0&-\frac{1}{2} &-\xi
\end{pmatrix}
= \Phi J(\xi) \Phi^{-1}
\end{equation*}
with
$\Phi=\text{diag}(1,1,1,\frac{1}{2(1+\xi^2)})  \in \Aut{\mathfrak{u}(2)}.$
\end{remark}

\begin{remark}
\rm
$\mathfrak{u}(2)$ has no abelian complex structures since, for $J(\xi),$
$\frak{m} =\Cmath \tilde{J}_1 \oplus \Cmath \tilde{J}_3$ is the solvable
Lie algebra
$[\tilde{J}_1 ,\tilde{J}_3]= i(1-i\xi)\tilde{J}_1.$
\end{remark}

\begin{corollary}
${\frak{X}}_{\mathfrak{u}(2)}$
is comprised of the matrices
\begin{equation}
\label{Xu2}
\begin{pmatrix}
(a^1_4)^2 c^2 \xi
&
(a^3_4+a^2_4a^1_4c\xi)c&
 (a^3_4 a^1_4c\xi -a^2_4)c& a^1_4\\
- (a^3_4- a^2_4 a^1_4 c \xi)c&
(a^2_4)^2 c^2 \xi
& (a^3_4 a^2_4 c \xi + a^1_4)c & a^2_4\\
 (a^3_4 a^1_4 c \xi + a^2_4)c &
 (a^3_4 a^2_4 c \xi - a^1_4)c &
(a^3_4)^2 c^2 \xi &
 a^3_4\\
-(\xi^2 +1)c^2 a^1_4 &
-(\xi^2 +1)c^2 a^2_4  &
-(\xi^2 +1)c^2 a^3_4   &
-\xi
\end{pmatrix}
\end{equation}
with the conditions
\begin{equation}
\label{Xu2conditions}
\xi \in \Rmath, \quad
\left(\begin{smallmatrix}
a^1_4\\a^2_4\\a^3_4
\end{smallmatrix}\right) \in \Rmath^3 \setminus \{0\} ,\quad
c= \pm \left(
(a^1_4)^2 +(a^2_4)^2 +(a^3_4)^2
\right)^{-\frac{1}{2}}.
\end{equation}
\end{corollary}
\begin{proof}
As is known, any $R \in SO(3)$ can be written
\begin{equation}
\label{R}
\begin{pmatrix}
u^2-v^2 -w^2+s^2 &-2(uv +ws)& 2(-uw +sv)\\
2(-sw+uv)&u^2-v^2 +w^2-s^2 &-2(su +vw)\\
2(sv+uw)&2(su-vw)&u^2+v^2 -w^2-s^2
\end{pmatrix}
\end{equation}
for $q=(u,v,w,s)\in \Smath^3$ ($R$ can be written in exactly 2 ways by means of
$q$ and $-q$).
Hence any $\Phi \in \Aut{\mathfrak{u}(2)}$ can be written
$$\Phi = \begin{pmatrix} &R&& \begin{matrix} 0\\0\\0 \end{matrix}\\ 0&0&0&c\end{pmatrix}$$
with $R$ as in (\ref{R}) and $c\in \Rmath^*.$
Then we get
for $\Phi J(\xi) \Phi^{-1}$
the matrix
(\ref{Xu2})
with
\begin{eqnarray}
\label{XU21}
a^1_4 &=& \frac{2}{c} (sv -uw),\\
\label{XU22}
a^2_4 &=& -\frac{2}{c} (su +vw),\\
\label{XU23}
a^3_4 &=& \frac{1}{c} (2u^2 +2v^2 -1).
\end{eqnarray}
From $u^2+v^2+w^2+s^2=1,$ one gets
$(a^1_4)^2 +(a^2_4)^2 +(a^3_4)^2 = \frac{1}{c^2}.$
Conversely,  for any matrix $J$ of the form
(\ref{Xu2})
with conditions
(\ref{Xu2conditions})
there exist
$\Phi \in \Aut{\mathfrak{u}(2)}$ and $\xi \in \Rmath$ such that
 $J= \Phi J(\xi) \Phi^{-1}.$  This amounts to the existence of
$q=(u,v,w,s)\in \Smath^3$ such that
equations
(\ref{XU21}),
(\ref{XU22}),
(\ref{XU23}) hold true,
and follows from the fact that the map
$\Smath^3 \rightarrow  \Smath^2$
$q \mapsto (c a^1_4,c a^2_4,c a^3_4)$  is the Hopf fibration.
\end{proof}

\begin{corollary}
${\frak{X}}_{\mathfrak{u}(2)}$
is a closed 4-dimensional (smooth) submanifold of $\Rmath^{16}$ with 2 connected components,
each of them diffeomorphic to
$\Rmath \times \left(\Rmath^3 \setminus \{0\}\right).$
\end{corollary}
\begin{proof}
Denote ${\frak{X}}_{\mathfrak{u}(2)}^+$
(resp. ${\frak{X}}_{\mathfrak{u}(2)}^-$) the subset of
those $ J \in {\frak{X}}_{\mathfrak{u}(2)}$  with $c>0$ (resp. $c<0$).
As $c$ is uniquely defined by  the matrix $J=(a^i_j) \in {\frak{X}}_{\mathfrak{u}(2)}$  by the formula
$$2c = \frac{a^3_4(a^1_2-a^2_1)+a^2_4(-a^1_3+a^3_1)+a^1_4(a^2_3-a^3_2)}{(a^1_4)^2 +(a^2_4)^2 +(a^3_4)^2 },$$
one has ${\frak{X}}_{\mathfrak{u}(2)} ={\frak{X}}_{\mathfrak{u}(2)}^+ \cup    {\frak{X}}_{\mathfrak{u}(2)}^-$ with disjoint union.
       ${\frak{X}}_{\mathfrak{u}(2)}^+$
(resp. ${\frak{X}}_{\mathfrak{u}(2)}^-$) is a closed subset of
$\Rmath^{16}.$
It hence suffices to prove that
${\frak{X}}_{\mathfrak{u}(2)}^+$
is a regular submanifold, the case of
${\frak{X}}_{\mathfrak{u}(2)}^-$ being analogous.
Let $F \, : \Rmath \times \left(\Rmath^3 \setminus \{0\}\right)  \rightarrow   {\frak{X}}_{\mathfrak{u}(2)}^+$ be the bijection
defined by $F(\xi,(a^1_4,a^2_4,a^3_4)) = J$ where $J$ is the matrix
(\ref{Xu2}) with
$c=  \left(
(a^1_4)^2 +(a^2_4)^2 +(a^3_4)^2
\right)^{-\frac{1}{2}}.$
We equip
${\frak{X}}_{\mathfrak{u}(2)}^+$
with the differentiable structure transferred from
$\Rmath \times \left(\Rmath^3 \setminus \{0\}\right).$
The injection $i$ from
${\frak{X}}_{\mathfrak{u}(2)}^+$ into the open subset $X \subset
\Rmath^{16}$ defined by
$(a^1_4)^2 +(a^2_4)^2 +(a^3_4)^2 \neq 0$ is smooth. Now, there is a smooth retraction
$r \, : \, X \mapsto
{\frak{X}}_{\mathfrak{u}(2)}^+$
defined by
$r(A) = F(-a^4_4, (a^1_4,a^2_4,a^3_4))$ for $A=(a^i_j) \in X.$ Hence $i$ is an immersion
and the topology of
${\frak{X}}_{\mathfrak{u}(2)}^+$ is the induced topology of
$\Rmath^{16}.$
\end{proof}
\section{$SU(2)\times SU(2)$.}
\begin{lemma}
\label{automsu22}
$\Aut{(\mathfrak{su}(2) \oplus \mathfrak{su}(2))} =
\left( SO(3) \times SO(3)\right) \cup \, \tau \left( SO(3) \times SO(3) \right)$
where
$\tau = \left(\begin{smallmatrix} 0&I& \\I&0 \end{smallmatrix}\right)$
is the switch between the two factors of
$\mathfrak{su}(2) \oplus \mathfrak{su}(2).$
\end{lemma}
\begin{proof}
Let $J_k^{(1)}$
$(1 \leqslant k \leqslant 3)$ (resp. $ J_\ell^{(2)}$
$(1 \leqslant \ell \leqslant 3))$
be the basis for the first (resp. the second) factor
 $\mathfrak{su}(2)^{(1)}$
 (resp. $\mathfrak{su}(2)^{(2)})$
of $\mathfrak{su}(2) \oplus \mathfrak{su}(2)$ with relations
(\ref{relationssu2}), and $\pi^{(1)}$ (resp. $\pi^{(2)}$) the corresponding projections.
Let $\Phi = \left(\begin{smallmatrix} \Phi_1&\Phi_2\\ \Phi_3&\Phi_4 \end{smallmatrix}\right) \in
\Aut{(\mathfrak{su}(2) \oplus \mathfrak{su}(2))} ,$
each $\Phi_j $ being a $3 \times 3$ matrix.
$\Phi_1 = \left(\pi^{(1)} \circ \Phi \right)_{|
 \mathfrak{su}(2)^{(1)} }$
is an homomorphism of
 $\mathfrak{su}(2)^{(1)}$ into itself. Hence the three columns of $\Phi_1$
 are two-by-two orthogonal vectors in $\Rmath^3$ and if one of them is zero, then the 3 of them are zero.
 In particular, if $\Phi_1 \neq 0$, then $\Phi_1 \in SO(3).$
 With the same reasoning, the same property holds true for  $\Phi_2, \Phi_3,\Phi_4.$
 Suppose first $\Phi_1 \neq 0.$
 For $k,\ell =1,2,3,$ $[\pi^{(1)}(\Phi(J^{(1)}_k)), \pi^{(1)}(\Phi(J^{(2)}_\ell))]
=  \pi^{(1)}(\Phi ([J^{(1)}_k, J^{(2)}_\ell])) =0.$
That implies that any column of $\Phi_1$ is colinear with any column of $\Phi_2,$ and hence
$\Phi_2=0$ since the columns of $\Phi_1$ are linearly independent.
Then $\det{\Phi_4} \neq 0,$ whence $\Phi_4 \in SO(3)$ and finally $\Phi_3=0$ by the above reasoning.
Hence $\Phi = \left(\begin{smallmatrix} \Phi_1&0\\ 0&\Phi_4 \end{smallmatrix}\right) \in
SO(3) \times SO(3).$
 Suppose now $\Phi_1 = 0.$
Then $\det{\Phi_2} \neq 0,$ whence $\Phi_2 \in SO(3),$
and $\det{\Phi_3} \neq 0,$ whence $\Phi_3 \in SO(3).$
By the same argument as before,
$\Phi_4=0.$
Hence $\Phi = \left(\begin{smallmatrix}0& \Phi_2\\ \Phi_3&0 \end{smallmatrix}\right)
=
\tau
\left(\begin{smallmatrix}\Phi_3&0\\ 0&\Phi_2\end{smallmatrix}\right)
\in\tau \left( SO(3) \times SO(3) \right).$
\end{proof}

\begin{theorem}
\label{torsionfreesu22}
Let  $J \, : \, {\mathfrak{su}(2) \oplus \mathfrak{su}(2)} \rightarrow
{\mathfrak{su}(2) \oplus \mathfrak{su}(2)}$ linear.
 $J$ has zero torsion, i.e. satisfies
(\ref{Jinv}), if and only if
there exists
$\Phi \in  SO(3) \times SO(3)$
such that
\begin{equation}
\label{Jsu2xsu2notorsion}
\Phi^{-1}J\Phi =
\begin{pmatrix}
0&1&0&0&0&0\\
-1&0&0&0&0&0\\
0&0&\xi^3_3&0&0&\xi^3_6\\
0&0&0&0&1&0\\
0&0&0&-1&0&0\\
0&0&\xi^6_3&0&0&\xi^6_6
\end{pmatrix}.
\end{equation}
\end{theorem}
\begin{proof}
From lemma \ref{product}  and   lemma \ref{su2notorsion},
there exists
$\Phi \in  SO(3) \times SO(3)$
such that
\begin{equation}
\label{17}
\Phi^{-1}J\Phi =
\begin{pmatrix}
0&1&0&\xi^1_4&\xi^1_5&\xi^1_6\\
-1&0&0&\xi^2_4&\xi^2_5&\xi^2_6\\
0&0&\xi^3_3&\xi^3_4&\xi^3_5&\xi^3_6\\
\xi^4_1&\xi^4_2&\xi^4_3&0&1&0\\
\xi^5_1&\xi^5_2&\xi^5_3&-1&0&0\\
\xi^6_1&\xi^6_2&\xi^6_3&0&0&\xi^6_6
\end{pmatrix}.
\end{equation}
Hence  we may suppose $J$ of the form (\ref{17}).
The matrix $(\xi^i_j)_{1\leqslant i \leqslant3,4\leqslant j \leqslant 6}$
(resp. $(\xi^i_j)_{4\leqslant i \leqslant6,1\leqslant j \leqslant 3}$)
is the matrix of
$\pi^1_2$
(resp. $\pi^2_1$)
of lemma \ref{product}.
Consider the vectors
$\bv{u} = {\pi}^1_2 J^{(2)}_1,
\bv{v} = {\pi}^1_2 J^{(2)}_2,
\bv{w} = {\pi}^1_2 J^{(2)}_3.$
From lemma \ref{su2notorsion} (ii) one has:
$$
[{\pi}^1_2 J^{(2)}_1,{\pi}^1_2 J^{(2)}_2] =
{\pi}^1_2 \, [J J^{(2)}_1,J^{(2)}_2] +
{\pi}^1_2 \, [J^{(2)}_1, J J^{(2)}_2]=
{\pi}^1_2 \, [-J^{(2)}_2,J^{(2)}_2] +
{\pi}^1_2 \, [J^{(2)}_1, J^{(2)}_1]=0,$$
 \begin{multline*}
[{\pi}^1_2 J^{(2)}_2,{\pi}^1_2 J^{(2)}_3] =
{\pi}^1_2 \, [J J^{(2)}_2,J^{(2)}_3] +
{\pi}^1_2 \, [J^{(2)}_2, J J^{(2)}_3]=
{\pi}^1_2 \, [J^{(2)}_1,J^{(2)}_3] +
{\pi}^1_2 \, [J^{(2)}_2, \xi^6_6 J^{(2)}_3]\\
=-{\pi}^1_2 \, J^{(2)}_2 +\xi^6_6
{\pi}^1_2 \, J^{(2)}_1,
 \end{multline*}

 \begin{multline*}
[{\pi}^1_2 J^{(2)}_1,{\pi}^1_2 J^{(2)}_3] =
{\pi}^1_2 \, [J J^{(2)}_1,J^{(2)}_3] +
{\pi}^1_2 \, [J^{(2)}_1, J J^{(2)}_3]=
{\pi}^1_2 \, [-J^{(2)}_2,J^{(2)}_3] +
{\pi}^1_2 \, [J^{(2)}_1, \xi^6_6 J^{(2)}_3]\\
=-{\pi}^1_2 \, J^{(2)}_1 -\xi^6_6  {\pi}^1_2 \, J^{(2)}_2.
 \end{multline*}
 That is:
\begin{eqnarray*}
\bv{u}\wedge \bv{v}&=&0,\\
\bv{v}\wedge \bv{w}&=& -\bv{v}+ \xi^6_6 \bv{u},\\
\bv{u}\wedge \bv{w}&=& -\bv{u}-\xi^6_6 \bv{v}
\end{eqnarray*}
which implies $\bv{u}=\bv{v}=0.$
With the same reasoning for $\pi^2_1,$ we get
\begin{equation}
J=
\begin{pmatrix}
0&1&0&0&0&\xi^1_6\\
-1&0&0&0&0&\xi^2_6\\
0&0&\xi^3_3&0&0&\xi^3_6\\
0&0&\xi^4_3&0&1&0\\
0&0&\xi^5_3&-1&0&0\\
0&0&\xi^6_3&0&0&\xi^6_6
\end{pmatrix}.
\end{equation}
Now      the  torsion equations
$16|3,26|3$
$36|4,36|5$ give the 2 Cramer systems
$\xi^2_6\xi^3_3 -\xi^1_6=0,$
$\xi^2_6 + \xi^3_3 \xi^1_6=0;$
$\xi^5_3\xi^6_6 +\xi^4_3=0,$
$-\xi^5_3 + \xi^6_6 \xi^4_3=0.$
Hence $\xi^1_6=\xi^2_6=\xi^4_3=\xi^5_3=0.$
Then all torsion equations vanish, and  the theorem  is proved
(\cite{companionarchive}, $\texttt{torsionsu22.red}$).
\end{proof}

\begin{corollary}
\label{equivCSsu22}
Any $J \in {\frak{X}}_{\mathfrak{su}(2) \oplus \mathfrak{su}(2)}$
is equivalent under some member of $SO(3)\times SO(3)$ to
\begin{equation}
\label{Jsu2xsu2}
J(\xi,\eta) =
\begin{pmatrix}
0&1&0&0&0&0\\
-1&0&0&0&0&0\\
0&0&\xi&0&0&\eta\\
0&0&0&0&1&0\\
0&0&0&-1&0&0\\
0&0&-\frac{1+{\xi}^2}{\eta}&0&0&-\xi
\end{pmatrix}
\end{equation}
with $\xi,\eta \in \Rmath, $ $\eta \neq 0.$
$J(\xi,\eta)$ and $J(\xi^{\prime},\eta^{\prime})$
are equivalent under some member of
$SO(3) \times SO(3)$
(resp. $\tau \left( SO(3) \times SO(3) \right)$)
if and only if
$\xi^{\prime}=\xi$
and $\eta^{\prime} =\eta$
(resp.
$\xi^{\prime}=-\xi$ and $\eta^{\prime} = -\frac{1+\xi^2}{\eta}$).
\end{corollary}
\begin{proof}
 $J$ in (\ref{Jsu2xsu2notorsion}) satisfies
 $J^2=-I$
 if and only if $\xi^3_6 \neq 0$ and
 $\xi^6_3 =-\frac{1+(\xi^3_3)^2}{\xi^3_6}, \xi^6_6=- \xi^3_3,$
leading to
$J(\xi,\eta)$
in (\ref{Jsu2xsu2})
with $\xi=\xi^3_3, \eta=\xi^3_6.$
\par
 Suppose
$J(\xi^{\prime},\eta^{\prime}) =
\Phi J(\xi,\eta)\Phi^{-1}$
with $\Phi =
\left(\begin{smallmatrix}\Phi_1&0\\ 0&\Phi_2\end{smallmatrix}\right)
\in
SO(3) \times SO(3).$
Then
$\left( \begin{smallmatrix}
0&1&0\\
-1&0&0\\
0&0&\xi^{\prime}
\end{smallmatrix} \right) =
\Phi_1
\left( \begin{smallmatrix}
0&1&0\\
-1&0&0\\
0&0&\xi
\end{smallmatrix} \right) \Phi_1^{-1}$
and
$\left( \begin{smallmatrix}
0&1&0\\
-1&0&0\\
0&0&-\xi^{\prime}
\end{smallmatrix} \right) =
\Phi_2
\left( \begin{smallmatrix}
0&1&0\\
-1&0&0\\
0&0&-\xi
\end{smallmatrix} \right) \Phi_2^{-1}
,$
which imply first $\xi^{\prime}=\xi$ and second
$\Phi_1=\text{diag}(R_1,1)$,
$\Phi_2=\text{diag}(R_2,1)$ with $R_1,R_2 \in SO(2).$
Then $
\left( \begin{smallmatrix}
0&0&0\\
0&0&0\\
0&0&\eta^{\prime}
\end{smallmatrix} \right) =
\Phi_1
\left( \begin{smallmatrix}
0&0&0\\
0&0&0\\
0&0&\eta
\end{smallmatrix} \right) \Phi_2^{-1} $
implies $\eta^{\prime} = \eta.$
\par
Now suppose
$J(\xi^{\prime},\eta^{\prime}) = \Psi J(\xi,\eta)\Psi^{-1}$
with $\Psi = \tau  \Phi
\in \tau \left(SO(3) \times SO(3)\right),$
$\Phi=
\left(\begin{smallmatrix}\Phi_1&0\\ 0&\Phi_2\end{smallmatrix}\right)
\in
SO(3) \times SO(3).$
Then
$\Phi J(\xi,\eta)\Phi^{-1} =
\tau J(\xi^{\prime},\eta^{\prime}) \tau = J(-\xi^{\prime},
-\frac{1+{\xi^{\prime}}^2}{\eta^{\prime}}$). Hence
$\xi =-\xi^{\prime}$ and
$\eta=-\frac{1+{\xi^{\prime}}^2}{\eta^{\prime}}$, i.e.
$\eta^{\prime}=-\frac{1+\xi^2}{\eta}$.
\end{proof}

\begin{remark}
\rm
Lemma 1 in \cite{daurtseva} states that a left invariant almost complex structure
on
${SU}(2) \times {SU}(2)$ is integrable if and only if it has the form
$A I_{a,c} A^{-1}$ with
$A \in SO(3) \times SO(3),$
$a \in \Rmath, c \in \Rmath^*$,
and
\begin{equation*}
I_{a,c}=
\begin{pmatrix}
\frac{a}{c} & 0 &0 & -\frac{a^2+c^2}{c} &0 &0\\
0 & 0 &-1 & 0 &0 &0\\
0 & 1 &0 & 0 &0 &0\\
\frac{1}{c} & 0 & 0 & -\frac{a}{c} &0 &0\\
0 & 0 & 0 & 0&0 &-1\\
0 & 0 &0  & 0 & 1&0
\end{pmatrix}.
\end{equation*}
One has $\Phi^{-1} I_{a,c} \Phi =J(\frac{a}{c},-\frac{a^2+c^2}{c})$ with
$\Phi = \text{diag}\left(
\left(
\begin{smallmatrix}
0 & 0 &-1 \\
0& 1 &0 \\
1 & 0 &0
\end{smallmatrix}
\right),
\left(
\begin{smallmatrix}
0 & 0 &-1 \\
0& 1 &0 \\
1 & 0 &0
\end{smallmatrix}
\right)
\right)
\in SO(3) \times SO(3).
$
\end{remark}

\begin{remark}
\rm
$\frak{su}(2) \times \frak{su}(2)$
has no abelian complex structures since, for $J(\xi,\eta),$
$\frak{m} =\Cmath \tilde{J}^{(1)}_1 \oplus \Cmath \tilde{J}^{(1)}_3 \oplus \Cmath \tilde{J}^{(2)}_1$ is the solvable
Lie algebra
$[\tilde{J}^{(1)}_1 ,\tilde{J}^{(1)}_3]= i(1-i\xi)\tilde{J}^{(1)}_1,$
$[\tilde{J}^{(1)}_3 ,\tilde{J}^{(2)}_1]= \frac{1+\xi^2}{\eta} \tilde{J}^{(2)}_1.$
\end{remark}

\begin{corollary}
\label{corollarysu22}
${\frak{X}}_
{\mathfrak{su}(2) \oplus \mathfrak{su}(2)}$
is comprised of the matrices
\begin{equation}
\label{Xsu2xsu2}
\begin{pmatrix}
\lambda_1^2\xi&
-\lambda_1 \mu_1\xi +\nu_1&
\lambda_1 \nu_1\xi +\mu_1&
\eta \lambda_1 \lambda_2&
-\eta \lambda_1 \mu_2&
\eta \lambda_1 \nu_2\\
-\lambda_1 \mu_1\xi -\nu_1&
\mu_1^2\xi&
\lambda_1 -\mu_1\nu_1\xi &
-\eta \mu_1 \lambda_2&
\eta \mu_1 \mu_2&
-\eta \mu_1 \nu_2\\
\lambda_1 \nu_1\xi -\mu_1&
-\lambda_1 -\mu_1\nu_1\xi &
\nu_1^2\xi&
\eta \nu_1\lambda_2 &
-\eta \nu_1\mu_2 &
\eta \nu_1 \nu_2\\
-\frac{\xi^2+1}{\eta}\lambda_1\lambda_2&
\frac{\xi^2+1}{\eta}\mu_1\lambda_2&
-\frac{\xi^2+1}{\eta}\nu_1\lambda_2&
-\lambda_2^2\xi&
\lambda_2 \mu_2\xi +\nu_2&
-\lambda_2 \nu_2\xi +\mu_2\\
\frac{\xi^2+1}{\eta}\lambda_1\mu_2&
-\frac{\xi^2+1}{\eta}\mu_1\mu_2&
\frac{\xi^2+1}{\eta}\nu_1\mu_2&
\lambda_2 \mu_2\xi -\nu_2&
-\mu_2^2\xi&
\lambda_2 + \mu_2\nu_2\xi\\
-\frac{\xi^2+1}{\eta}\lambda_1\nu_2&
\frac{\xi^2+1}{\eta}\mu_1\nu_2&
-\frac{\xi^2+1}{\eta}\nu_1\nu_2&
-\lambda_2 \nu_2\xi -\mu_2&
-\lambda_2 +\mu_2\nu_2\xi&
-\nu_2^2\xi
\end{pmatrix}
\end{equation}
with
\begin{equation}
\label{Xsu2xsu2conditions}
(\xi,\eta) \in   \Rmath \times \Rmath^*, \quad
\left(\begin{smallmatrix} \lambda_i\\ \mu_i \\  \nu_i \end{smallmatrix}\right) \in \Smath^2 \quad i=1,2.
\end{equation}
\end{corollary}
\begin{proof}
${\frak{X}}_{\mathfrak{su}(2) \oplus \mathfrak{su}(2)}$
is comprised of the matrices $\Phi J(\xi,\eta) \Phi^{-1}, \; (\xi,\eta) \in
\Rmath \times \Rmath^*, \; \Phi \in
SO(3) \times SO(3) .$
Let
$\Phi =    \left(
\begin{smallmatrix}
\Phi_1&0\\0&\Phi_2
\end{smallmatrix}\right) \in
SO(3) \times SO(3).$
$\Phi_1,\Phi_2$ can be written in the form
(\ref{R})
for resp. $q_1=(u_1,v_1,w_1,s_1),\, q_2=(u_2,v_2,w_2,s_2) \in \Smath^3.$
Then  $\Phi J(\xi, \eta) \Phi^{-1} $ is the matrix
(\ref{Xsu2xsu2})
with  for $i=1,2$
\begin{eqnarray}
\label{Xsu2xsu21}
\lambda_i &=& 2 (s_iv_i -u_iw_i),\\
\label{Xsu2xsu22}
\mu_i &=& 2(s_i u_i +v_iw_i),\\
\label{Xsu2xsu23}
\nu_i &=&  2u_i^2 +2v_i^2 -1.
\end{eqnarray}
One has $\lambda_i^2+ \mu_i^2+ \nu_i^2=1.$
Conversely,  for any matrix $J$ of the form
(\ref{Xsu2xsu2})
with condition
(\ref{Xsu2xsu2conditions})
there exist
$\Phi \in
SO(3) \times SO(3)$
and
$(\xi,\eta) \in   \Rmath \times \Rmath^*$
such that
 $J= \Phi J(\xi,\eta) \Phi^{-1}.$  This amounts to the existence
 for $i=1,2$
of $q_i=(u_i,v_i,w_i,s_i)\in \Smath^3$ such that
equations
(\ref{Xsu2xsu21}),
(\ref{Xsu2xsu22}),
(\ref{Xsu2xsu23})
hold true, which again follows from the Hopf fibration.
\end{proof}

\begin{corollary}
${\frak{X}}_
{\mathfrak{su}(2) \oplus \mathfrak{su}(2)}$
is a  closed 6-dimensional (smooth) submanifold of $\Rmath^{36}$
diffeomorphic to
$\Rmath \times \Rmath^*  \times \left(\Smath^2\right)^2.$
\end{corollary}
\begin{proof}
Let $X$  the open subset
of $\Rmath^{36}$ of those matrices
$\left( a^i_j \right) _{1\leqslant i,j \leqslant 6}$ such that
\linebreak[4]
$H^2N_1N_2 \neq 0$ where
$H^2=\sum_{i=1}^{3} \sum_{j=4}^{6} \left( a^i_j \right) ^2 ,$
$N_1=\left(a^2_3-a^3_2\right)^2
+\left(a^1_3-a^3_1\right)^2
+\left(a^1_2-a^2_1\right)^2 ,$
$N_2=\left(a^5_6-a^6_5\right)^2
+\left(a^4_6-a^6_4\right)^2
+\left(a^4_5-a^5_4\right)^2$
and consider
 $F \, : \Rmath \times \Rmath^*  \times
\left(\Smath^2\right)^2  \rightarrow X$
defined by
$F(\xi,\eta,
(\lambda_1,\mu_1,\nu_1),
(\lambda_2,\mu_2,\nu_2)
) = J$ where $J$ is the matrix
(\ref{Xsu2xsu2}).

Observe first that $F$ is injective. In fact,
$\xi,
(\lambda_1,\mu_1,\nu_1),
(\lambda_2,\mu_2,\nu_2)$ can be retrieved from
$(a^i_j)=
F(\xi,\eta,
(\lambda_1,\mu_1,\nu_1),
(\lambda_2,\mu_2,\nu_2)
)$
by the formulas
$\xi = a^1_1+a^2_2+a^3_3,$ \,
\linebreak[4]
$(\lambda_1,\mu_1,\nu_1)
=
\left(
\frac{a^2_3-a^3_2}{\sqrt{N_1}},
\frac{a^1_3-a^3_1}{\sqrt{N_1}},
\frac{a^1_2-a^2_1}{\sqrt{N_1}}
\right) $,
$(\lambda_2,\mu_2,\nu_2)=
\left(
\frac{a^5_6-a^6_5}{\sqrt{N_2}},
\frac{a^4_6-a^6_4}{\sqrt{N_2}},
\frac{a^4_5-a^5_4}{\sqrt{N_2}}
\right);$
hence
\linebreak[4]
$
F(\xi,\eta,
(\lambda_1,\mu_1,\nu_1),
(\lambda_2,\mu_2,\nu_2)
) =
F(\xi^{\prime},\eta^{\prime},
(\lambda_1^{\prime},\mu_1^{\prime},\nu_1^{\prime}),
(\lambda_2^{\prime},\mu_2^{\prime},\nu_2^{\prime})
)$
implies $\xi=\xi^{\prime},
(\lambda_1^{\prime},\mu_1^{\prime},\nu_1^{\prime})
=(\lambda_1,\mu_1,\nu_1),
(\lambda_2^{\prime},\mu_2^{\prime},\nu_2^{\prime})
=
(\lambda_2,\mu_2,\nu_2),$
and then $\eta=\eta^{\prime}$ since
$\left(
\begin{smallmatrix}
 \lambda_1 \lambda_2&
- \lambda_1 \mu_2&
 \lambda_1 \nu_2\\
- \mu_1 \lambda_2&
 \mu_1 \mu_2&
- \mu_1 \nu_2\\
 \nu_1\lambda_2 &
- \nu_1\mu_2 &
 \nu_1 \nu_2
\end{smallmatrix}
\right)\neq 0.$
From the injectivity of $F$,
${\frak{X}}_
{\mathfrak{su}(2) \oplus \mathfrak{su}(2)} =
{\frak{X}}_
{\mathfrak{su}(2) \oplus \mathfrak{su}(2)}^+ \cup
{\frak{X}}_
{\mathfrak{su}(2) \oplus \mathfrak{su}(2)}^-$
 with disjoint union,
where
${\frak{X}}_
{\mathfrak{su}(2) \oplus \mathfrak{su}(2)}^\epsilon$
denotes the set of those  $J$s
having $\eta$ the sign of $\epsilon$ ($\epsilon=\pm$).
Now, the map $G_\epsilon  : X \rightarrow
 \Rmath \times \Rmath_\epsilon^*  \times
\left(\Smath^2\right)^2$
defined by
\begin{multline}
G_\epsilon(\left( a^i_j \right) ) =
\left(
a^1_1+a^2_2+a^3_3, \,
\epsilon \sqrt{\sum_{i=1}^{3} \sum_{j=4}^{6} \left( a^i_j \right) ^2} ,\,
\left(
\frac{a^2_3-a^3_2}{\sqrt{N_1}},
\frac{a^1_3-a^3_1}{\sqrt{N_1}},
\frac{a^1_2-a^2_1}{\sqrt{N_1}}
\right) \right.,
\\
\left.
\left(
\frac{a^5_6-a^6_5}{\sqrt{N_2}},
\frac{a^4_6-a^6_4}{\sqrt{N_2}},
\frac{a^4_5-a^5_4}{\sqrt{N_2}}
\right)
\right)
\end{multline}
is a smooth retraction for the restriction
$F_\epsilon$ of $F$ to
 $\Rmath \times \Rmath_\epsilon^*  \times
\left(\Smath^2\right)^2.$
Hence
$F_\epsilon$ is an immersion and the topology of
${\frak{X}}_{\mathfrak{su}(2) \oplus \mathfrak{su}(2)}^\epsilon$
is the induced topology from $X.$
The corollary follows.
\end{proof}

\begin{remark}
\rm
We may consider
$\mathfrak{u}(2)$
as a  subalgebra
of ${\mathfrak{su}(2) \oplus \mathfrak{su}(2)}$
by identifying $J_1,J_2,J_3,J_4$ to $J^{(1)}_1,
J^{(1)}_2,
J^{(1)}_3,
J^{(2)}_3$ respectively.
Then the complex structure $J$ in
(\ref{Xsu2xsu2}) leaves
$\mathfrak{u}(2)$ invariant if and only if $\lambda_2=\mu_2=0, \nu_2=\pm1.$
For the  restriction of $J$  to
$\mathfrak{u}(2)$  to be
(\ref{Xu2}), one must take
$\lambda_1=\frac{a^1_4}{\eta\nu_2},
\mu_1=-\frac{a^2_4}{\eta\nu_2},
\nu_1=\frac{a^3_4}{\eta\nu_2}$
with $c=\frac{\nu_2}{\eta}.$
Then
\begin{equation*}
J=
\begin{pmatrix}
(a^1_4)^2 c^2 \xi
&
(a^3_4+a^2_4a^1_4c\xi)c&
 (a^3_4 a^1_4c\xi -a^2_4)c& 0&0&a^1_4\\
- (a^3_4- a^2_4 a^1_4 c \xi)c&
(a^2_4)^2 c^2 \xi
& (a^3_4 a^2_4 c \xi + a^1_4)c &0&0& a^2_4\\
 (a^3_4 a^1_4 c \xi + a^2_4)c &
 (a^3_4 a^2_4 c \xi - a^1_4)c &
(a^3_4)^2 c^2 \xi &0&0&
 a^3_4\\
 0&0&0&0&\nu_2&0\\
 0&0&0&-\nu_2&0&0\\
-(\xi^2 +1)c^2 a^1_4 &
-(\xi^2 +1)c^2 a^2_4  &
-(\xi^2 +1)c^2 a^3_4   &
0&0&
-\xi
\end{pmatrix}.
\end{equation*}
Hence any complex structure on
$\mathfrak{u}(2)$
can be extended in 2 (in general non equivalent) ways to a complex structure on
 ${\mathfrak{su}(2) \oplus \mathfrak{su}(2)}.$
 For example, $J(\xi)$ can be extended
 (here $a^1_4=a^2_4=0,a^3_4=1,c=1$)
 with $\nu_2=1$
 to $J(\xi,1)$
 or
 with $\nu_2=-1$
 to
$\left(\begin{smallmatrix}
0&1&0&0&0&0\\
-1&0&0&0&0&0\\
0&0&\xi&0&0&1\\
0&0&0&0&-1&0\\
0&0&0&1&0&0\\
0&0&-(1+\xi^2)&0&0&-\xi
\end{smallmatrix} \right)$
which is equivalent to $J(\xi,-1).$
Now,
$J(\xi,-1) \cong J(\xi,1)
\Leftrightarrow  \xi=0.$
\end{remark}

\section{$SU(2)^N.$}
The results of
lemma \ref{automsu22},
theorem \ref{torsionfreesu22}
and corollary \ref{equivCSsu22}
easily
generalize
in the following way.
\begin{lemma}
For any $N \in \Nmath^*,$
$\Aut{(\mathfrak{su}(2))^N} =
SO(3)^N \cup \, \left( \bigcup_{\sigma \in \Sigma} \, \tau_\sigma \, \left( SO(3)^N \right) \right)$
(disjoint reunion)
where
$\Sigma$ is the set of circular permutations of $\{1, \cdots, N\}$ having no fixed point,
and $\tau_\sigma=(T^i_j)_{1\leqslant i,j \leqslant N}$
with the $T^i_j$s the $3\times 3$ blocks
$T^i_j= \delta_{i,\sigma(j)}\, I$ ($I$ the $3\times 3$ identity and $\delta_{k,\ell}$ the Kronecker symbol).
\end{lemma}

\begin{theorem}
Let  $J \, : \, \mathfrak{su}(2)^N  \rightarrow
\mathfrak{su}(2)^N$ linear.
 $J$ has zero torsion if and only if
there exists
$\Phi \in  SO(3)^N$
and $M=(\xi^{3i}_{3j})_{1\leqslant i,j \leqslant N} \in \text{gl}(N,\Rmath)$
such that
$\Phi^{-1}J\Phi = J(M)$
with $J(M) =(J^i_j(M))_{1\leqslant i,j \leqslant N}$
and the $J^i_j(M)$s the following $3\times 3$ blocks
\begin{equation}
\label{torsionfreesuN}
J^i_i(M) =
\begin{pmatrix}
0&1&0\\
-1&0&0\\
0&0&\xi^{3i}_{3i}
\end{pmatrix} \quad (1\leqslant i \leqslant N)\; , \quad
J^i_j(M) =
\begin{pmatrix}
0&0&0\\
0&0&0\\
0&0&\xi^{3i}_{3j}
\end{pmatrix} \quad (1\leqslant i,j \leqslant N, \; i \neq j).
\end{equation}
\end{theorem}
(Here we used that the analogs of
$16|3,26|3$
$36|4,36|5$
at the end of the proof of Theorem
\ref{torsionfreesu22}
are resp. , with $i <j,$
\begin{eqnarray*}
3i-2,3j|3i : &\xi^{3i-2}_{3j} -\xi^{3i}_{3i}\xi^{3i-1}_{3j}=0\\
3i-1,3j|3i : &\xi^{3i-1}_{3j} +\xi^{3i}_{3i}\xi^{3i-2}_{3j}=0\\
3i,3j|3j-2 : &-\xi^{3j-2}_{3i} -\xi^{3j}_{3j}\xi^{3j-1}_{3i}=0\\
3i-1,3j|3j-1 : &\xi^{3j-1}_{3i} +\xi^{3j}_{3j}\xi^{3j-2}_{3i}=0
\end{eqnarray*}
and give
$\xi^{3i-2}_{3j}=\xi^{3i-1}_{3j}=
\xi^{3j-2}_{3i}=\xi^{3j-1}_{3i}=0.$ Then all torsion equations vanish.)

\begin{example}
\rm
For $N=4,$
$$J(M)=
\left(
\begin{array}{ccc|ccc||ccc|ccc}
 0&1&0&0&0&0& 0&0&0&0&0&0\\
 -1&0&0&0&0&0& 0&0&0&0&0&0\\
0&0&\boxed{\xi^3_3}&0&0&\boxed{\xi^3_6} &0&0&\boxed{\xi^3_9}&0&0&\boxed{\xi^3_{12}}\\
\hline
0&0&0&0&1&0 & 0&0&0&0&0&0\\
0&0&0&-1&0&0 & 0&0&0&0&0&0\\
0&0&\boxed{\xi^6_3}&0&0&\boxed{\xi^6_6}&0&0&\boxed{\xi^6_9}&0&0&\boxed{\xi^6_{12}}\\
\hline
\hline
 0&0&0&0&0&0 &0&1&0&0&0&0\\
 0&0&0&0&0&0 &-1&0&0&0&0&0\\
0&0&\boxed{\xi^9_3}&0&0&\boxed{\xi^9_{6}}&0&0&\boxed{\xi^9_9}&0&0&\boxed{\xi^9_{12}}\\
\hline
 0&0&0&0&0&0& 0&0&0&0&1&0\\
 0&0&0&0&0&0& 0&0&0&-1&0&0\\
0&0&\boxed{\xi^{12}_3}&0&0&\boxed{\xi^{12}_{6}}&0&0&\boxed{\xi^{12}_9}&0&0&\boxed{\xi^{12}_{12}}

\end{array}
\right)
$$
$$M=
\begin{pmatrix}
\xi^3_3&\xi^3_6&\xi^3_9&\xi^3_{12}\\
\xi^6_3&\xi^6_6&\xi^6_9&\xi^6_{12}\\
\xi^9_3&\xi^9_6&\xi^9_9&\xi^9_{12}\\
\xi^{12}_3&\xi^{12}_6&\xi^{12}_9&\xi^{12}_{12}
\end{pmatrix}
.$$
\end{example}
\begin{corollary}
For even $N,$ any $J \in {\frak{X}}_{\mathfrak{su}(2)^N}$
is equivalent
under some member of $SO(3)^N$
to some
$J(M) = (J^i_j(M))_{1\leqslant i,j \leqslant N}$
with
$M=\left(\xi^{3i}_{3j}\right)_{1\leqslant i,j \leqslant N}$
such that $M^2=-I$
and
$J^i_j(M)$ defined in
(\ref{torsionfreesuN}).
$J(M)$ and $J(M^{\prime})$
are equivalent under some member of
$SO(3)^N$
(resp. $\tau_\sigma \left( SO(3)^N \right), \, \sigma \in \Sigma$)
if and only if
$M^{\prime}=M$
(resp.
$M^{\prime}=M^{\sigma^{-1}},$  $M^{\sigma^{-1}} = \left(\xi^{3\,\sigma^{-1}(i)}_{3\,\sigma^{-1}(j)}\right)_{1\leqslant i,j \leqslant N}$ (here we make use of $(\tau_\sigma)^{-1}= \tau_{\sigma^{-1}}$ and
$\tau_\sigma J(M^{\prime}) (\tau_\sigma)^{-1}=\left(J^{\sigma^{-1}(i)}_{\sigma^{-1}(j)}(M
^{\prime}) \right)_{1\leqslant i,j \leqslant N}$
).
\end{corollary}
\begin{example}
For $N=2,$ $\Sigma$ consists only of the transposition $(1,2);$
$M=\begin{pmatrix} \xi^3_3 &\xi^3_6\\ -\frac{1+(\xi^3_3)^2}{\xi^3_6}&-\xi^3_3
\end{pmatrix},$
$M^{\prime}=\begin{pmatrix} {{\xi^{\prime}}^3_3} &{{\xi^{\prime}}^3_6}\\
-\frac{1+({{\xi^{\prime}}^3_3})^2}{{{\xi^{\prime}}^3_6}}&-{{\xi^{\prime}}^3_3}
\end{pmatrix}.$
For $\sigma=(1,2),$
the condition
$M^{\prime}=M^{\sigma^{-1}}$ reads
${{\xi^{\prime}}^3_3} = -\xi^3_3,
\,
{{\xi^{\prime}}^3_6} =
-\frac{1+({\xi^3_3})^2}{{\xi^3_6}}$
and is that of  Corollary \ref{equivCSsu22}.
\end{example}
\section{$U(2)\times U(2)$.}
\begin{lemma}
\label{automu22}
$\Aut{(\mathfrak{u}(2) \oplus \mathfrak{u}(2))} =
H \cup \, \tau H$
where
$\tau = \left(\begin{array}{c|c} 0&I \\ \hline I&0 \end{array}\right)$
is the switch between the two factors of
$\mathfrak{u}(2) \oplus \mathfrak{u}(2),$
$H=
\left\{
\left(
\begin{array}{cc|cc}
\Phi_1&0&0&0 \\
0&b^4_4&0& b^4_8\\
\hline
0&0 & \Phi_4&0 \\
0&b^8_4 &0 &b^8_8
\end{array}
\right)
, \Phi_1,\Phi_4 \in SO(3), b^4_4b^8_8-b^4_8b^8_4 \neq 0
\right\},$
$\tau H=
\left\{
\left(
\begin{array}{cc|cc}
0&0&\Phi_2&0 \\
0&b^4_4&0& b^4_8\\
\hline
\Phi_3&0 & 0&0 \\
0&b^8_4 &0 &b^8_8
\end{array}
\right)
, \Phi_2,\Phi_3 \in SO(3), b^4_4b^8_8-b^4_8b^8_4 \neq 0
\right\}.$
\end{lemma}
\begin{proof}
Analog to that of lemma \ref{automsu22}.
\end{proof}
\begin{theorem}
\label{torsionfreeu22}
(i) Let  $J \, : \,
{\mathfrak{u}(2) \oplus \mathfrak{u}(2)}
\rightarrow
{\mathfrak{u}(2) \oplus \mathfrak{u}(2)}
$ linear.
 $J$ has zero torsion
if and only if
there exists
$\Phi \in  \left( SO(3)\times \Rmath_+^*\right)^2 \subset H$
and $M \in gl(4,\Rmath)$
such that
$\Phi^{-1}J\Phi = K(M),$ where
\begin{equation}
\label{Ju2xu2notorsion}
K(M) =
\left(
\begin{array}{cccc|cccc}
0&1&0&0&0&0&0&0\\
-1&0&0&0&0&0&0&0\\
0&0&\xi^3_3&\xi^3_4&0&0&\xi^3_7&\xi^3_8\\
0&0&\xi^4_3&\xi^4_4&0&0&\xi^4_7&\xi^4_8\\
\hline
0&0&0&0&0&1&0&0\\
0&0&0&0&-1&0&0&0\\
0&0&\xi^7_3&\xi^7_4&0&0&\xi^7_7&\xi^7_8\\
0&0&\xi^8_3&\xi^8_4&0&0&\xi^8_7&\xi^8_8
\end{array}\right) ,
\quad
M=
\begin{pmatrix}
\xi^3_3&\xi^3_4&\xi^3_7&\xi^3_8\\
\xi^4_3&\xi^4_4&\xi^4_7&\xi^4_8\\
\xi^7_3&\xi^7_4&\xi^7_7&\xi^7_8\\
\xi^8_3&\xi^8_4&\xi^8_7&\xi^8_8
\end{pmatrix}.
\end{equation}
\\ (ii)
For $M,M^{\prime}\in gl(4,\Rmath),$ there exists some
$\Phi  \in H$
such that $K(M^{\prime})=\Phi K(M)\Phi^{-1}$ if and only if
there exists
$\left(\begin{smallmatrix}b^4_4&b^4_8\\b^8_4&b^8_8 \end{smallmatrix}\right) \in GL(2,\Rmath)$
such that
$M^{\prime}= G MG^{-1},$ with
\begin{equation}
\label{G}
G=
\begin{pmatrix}
1&0&0&0\\
0&b^4_4&0&b^4_8\\
0&0&1&0\\
0&b^8_4&0&b^8_8
\end{pmatrix} \in GL(4,\Rmath).
\end{equation}
\\ (iii)
For $M,M^{\prime}\in gl(4,\Rmath),$ there exists
$\Psi  \in \tau H$
such that $K(M^{\prime})=\Psi K(M)\Psi^{-1}$ if and only if
there exists
$\Phi  \in H$
such that $K(M^{\prime})=\Phi K(M)\Phi^{-1}.$
\end{theorem}
\begin{proof}
(i)
From lemma \ref{product}  and
theorem \ref{theoremu2notorsion}(i),
there exists
$\Phi \in  \left( SO(3)\times \Rmath_+^*\right)^2 \subset H$
such that
\begin{equation}
\label{29}
\Phi^{-1}J\Phi =
\left(
\begin{array}{cccc|cccc}
0&1&0&0&\xi^1_5&\xi^1_6&\xi^1_7&\xi^1_8\\
-1&0&0&0&\xi^2_5&\xi^2_6&\xi^2_7&\xi^2_8\\
0&0&\xi^3_3&\xi^3_4&\xi^3_5&\xi^3_6&\xi^3_7&\xi^3_8\\
0&0&\xi^4_3&\xi^4_4&\xi^3_5&\xi^3_6&\xi^4_7&\xi^4_8\\
\hline
\xi^5_1&\xi^5_2&\xi^5_3&\xi^5_4&0&1&0&0\\
\xi^6_1&\xi^6_2&\xi^6_3&\xi^6_4&-1&0&0&0\\
\xi^7_1&\xi^7_2&\xi^7_3&\xi^7_4&0&0&\xi^7_7&\xi^7_8\\
\xi^8_1&\xi^8_2&\xi^8_3&\xi^8_4&0&0&\xi^8_7&\xi^8_8
\end{array}\right).
\end{equation}
Hence we may suppose $J$ of the form (\ref{29}).
The matrix $(\xi^i_j)_{1\leqslant i \leqslant4,5\leqslant j \leqslant 8}$
(resp. $(\xi^i_j)_{5\leqslant i \leqslant8,1\leqslant j \leqslant 4}$)
is the matrix of
$\pi^1_2$
(resp. $\pi^2_1$)
of lemma \ref{product}.
Consider
$\bv{u} = {\pi}^1_2 J^{(2)}_1,
\bv{v} = {\pi}^1_2 J^{(2)}_2,
\bv{w} = {\pi}^1_2 J^{(2)}_3,
\bv{z} = {\pi}^1_2 J^{(2)}_4.$
From lemma \ref{su2notorsion} (ii) one has:
\begin{eqnarray*}
\left[ \bv{u}, \bv{v}\right] &=&0,\\
\left[ \bv{v}, \bv{w}\right] &=& -\bv{v}+ \xi^7_7 \bv{u},\\
\left[ \bv{u}, \bv{w}\right] &=& -\bv{u}-\xi^7_7 \bv{v},\\
\left[ \bv{u}, \bv{z}\right] &=& -\xi^7_8 \bv{v}\\
\left[ \bv{v}, \bv{z}\right] &=& \xi^7_8 \bv{u}\\
\left[ \bv{w}, \bv{z}\right] &=& 0
\end{eqnarray*}
which implies $\bv{u}=\bv{v}=0.$
With the same reasoning for $\pi^2_1,$ we get
\begin{equation*}
J=
\left(
\begin{array}{cccc|cccc}
0&1&0&0&0&0&\xi^1_7&\xi^1_8\\
-1&0&0&0&0&0&\xi^2_7&\xi^2_8\\
0&0&\xi^3_3&\xi^3_4&0&0&\xi^3_7&\xi^3_8\\
0&0&\xi^4_3&\xi^4_4&0&0&\xi^4_7&\xi^4_8\\
\hline
0&0&\xi^5_3&\xi^5_4&0&1&0&0\\
0&0&\xi^6_3&\xi^6_4&-1&0&0&0\\
0&0&\xi^7_3&\xi^7_4&0&0&\xi^7_7&\xi^7_8\\
0&0&\xi^8_3&\xi^8_4&0&0&\xi^8_7&\xi^8_8
\end{array}\right).
\end{equation*}
Now      the  torsion equations
$17|3,27|3,$
$18|3,28|3,$
$35|7,36|7,$
$45|7,46|7,$
give the 4 Cramer systems
$\xi^2_7\xi^3_3 -\xi^1_7=0,$
$\xi^2_7 + \xi^3_3 \xi^1_7=0;$
$\xi^2_8\xi^3_3 -\xi^1_8=0,$
$\xi^2_8 + \xi^3_3 \xi^1_8=0;$
$\xi^6_3\xi^7_7 -\xi^5_3=0,$
$\xi^6_3 +\xi^7_7\xi^5_3=0,$
$\xi^6_4\xi^7_7 -\xi^5_4=0,$
$\xi^6_4 +\xi^7_7\xi^5_4=0.$

Hence $\xi^1_7=\xi^2_7=\xi^1_8=\xi^2_8=0=
\xi^5_3=\xi^6_3=\xi^5_4=\xi^6_4=0.$
Then all torsion equations vanish
(\cite{companionarchive} $\texttt{torsionu22.red}$).
\\(ii)
Suppose there exists
$$ \Phi =\left(
\begin{array}{cc|cc}
\Phi_1&0&0&0 \\
0&b^4_4&0& b^4_8\\
\hline
0&0 & \Phi_4&0 \\
0&b^8_4 &0 &b^8_8
\end{array}
\right) \in H
,\quad  \Phi_1,\Phi_4 \in SO(3), b^4_4b^8_8-b^4_8b^8_4 \neq 0
$$
such that $\Phi K(M) = K(M^{\prime})\Phi,$
with
$
M=\left(
\begin{smallmatrix}
\xi^3_3&\xi^3_4&\xi^3_7&\xi^3_8\\
\xi^4_3&\xi^4_4&\xi^4_7&\xi^4_8\\
\xi^7_3&\xi^7_4&\xi^7_7&\xi^7_8\\
\xi^8_3&\xi^8_4&\xi^8_7&\xi^8_8
\end{smallmatrix}\right),$
$M^{\prime}=\left(
\begin{smallmatrix}
{\xi^{\prime}}^3_3  &{\xi^{\prime}}^3_4  &{\xi^{\prime}}^3_7   &{\xi^{\prime}}^3_8\\
{\xi^{\prime}}^4_3  &{\xi^{\prime}}^4_4  &{\xi^{\prime}}^4_7   &{\xi^{\prime}}^4_8\\
{\xi^{\prime}}^7_3  &{\xi^{\prime}}^7_4  &{\xi^{\prime}}^7_7   &{\xi^{\prime}}^7_8\\
{\xi^{\prime}}^8_3  &{\xi^{\prime}}^8_4  &{\xi^{\prime}}^8_7   &{\xi^{\prime}}^8_8
\end{smallmatrix}\right)$.
One has
$\Phi K(M) = \left(\begin{smallmatrix}A&B\\C&D\end{smallmatrix}\right)$
with
$A = \left(\begin{array}{ccc|c}
&A_1&&*\\
\hline
*&*&*&*
\end{array}
\right),$
$D = \left(\begin{array}{ccc|c}
&D_1&&*\\
\hline
*&*&*&*
\end{array}
\right),$
$A_1 = \Phi_1\left(\begin{smallmatrix}0&1&0\\-1&0&0\\0&0&\xi^3_3\end{smallmatrix}\right),$
$D_1 = \Phi_4\left(\begin{smallmatrix}0&1&0\\-1&0&0\\0&0&\xi^7_7\end{smallmatrix}\right),$
and $K(M^{\prime})\Phi = \left(\begin{smallmatrix}A^{\prime}&B^{\prime}\\C^{\prime}&D^{\prime}\end{smallmatrix}\right)$
with
$A^{\prime} = \left(\begin{array}{ccc|c}
&A^{\prime}_1&&*\\
\hline
*&*&*&*
\end{array}
\right),$
$D^{\prime} = \left(\begin{array}{ccc|c}
&D^{\prime}_1&&*\\
\hline
*&*&*&*
\end{array}
\right),$
$A^{\prime}_1 = \left(\begin{smallmatrix}0&1&0\\-1&0&0\\0&0&{\xi^{\prime}}^3_3\end{smallmatrix}\right)\Phi_1,$
$D^{\prime}_1 = \left(\begin{smallmatrix}0&1&0\\-1&0&0\\0&0&{\xi^{\prime}}^7_7\end{smallmatrix}\right)\Phi_4.$
Hence
${\xi^{\prime}}^3_3=\xi^3_3,$
${\xi^{\prime}}^7_7=\xi^7_7,$ and
$\Phi_1= \text{diag}(R_1,1),$
$\Phi_4= \text{diag}(R_4,1),$ $R_1,R_4 \in SO(2).$
Now
$$ \Phi =\left(
\begin{array}{cccccc}
R_1& 0& 0 &0 &0 &0 \\
0&1&0&0&0&0\\
  0  & 0&b^4_4&0 &0  &b^4_8 \\
  0  &0&0 & R_4&0 & 0\\
   0 &0&0 &0 &1 & 0\\
   0 &0&b^8_4 &0 &0 &b^8_8
\end{array}
\right),$$
and we can reindex the basis as the new basis
$(J^{(1)}_1,
J^{(1)}_2,
J^{(2)}_1,
J^{(2)}_2,
J^{(1)}_3,
J^{(1)}_4,
J^{(2)}_3,
J^{(2)}_4),$ so that $\Phi$, $K(M),$ $K(M^{\prime})$ have respectives matrices
in the new basis:
\begin{multline*}
\Phi =\left(
\begin{array}{cc|cccc}
R_1& 0& 0 &0 &0 &0 \\
0&R_4& 0 &0 &0 &0 \\
\hline
0&0& 1 &0 &0 &0 \\
 0&0& 0  & b^4_4&0 &b^4_8 \\
   0 &0 &0 &0 & 1& 0\\
 0&0& 0  & b^8_4&0 &b^8_8
\end{array}
\right), \;
K(M) =\left(
\begin{array}{cccc|c}
0& 1& 0 &0 &0 \\
-1& 0 &0 &0 &0 \\
0&0& 0 &1 &0 \\
0&0& -1 &0 &0 \\
\hline
0&0& 0 &0 &M
\end{array}
\right),\\
K(M^{\prime}) =\left(
\begin{array}{cccc|c}
0& 1& 0 &0 &0 \\
-1& 0 &0 &0 &0 \\
0&0& 0 &1 &0 \\
0&0& -1 &0 &0 \\
\hline
0&0& 0 &0 &M^{\prime}
\end{array}
\right).
\end{multline*}
The conclusion follows.
\\(iii)
One has $\tau K(M^{\prime}) \tau = K(\tau_1 M^{\prime} \tau_1)$
with
$\tau_1$ the same as $\tau$ yet with $2\times 2$ blocks:
$\tau_1=
\left(\begin{smallmatrix}
0&0&1&0\\
0&0&0&1\\
1&0&0&0\\
0&1&0&0\\
\end{smallmatrix}\right).$
Now, let
$\Psi= \tau \Phi  \in \tau H,$ $\Phi\in H.$
Then
$K(M^{\prime})=\Psi K(M)\Psi^{-1}$ if and only if
$\Phi K(M)\Phi^{-1} = \tau K(M^{\prime}) \tau=
K(\tau_1 M^{\prime} \tau_1)$ i.e. there exists
$\left(\begin{smallmatrix}b^4_4&b^4_8\\b^8_4&b^8_8 \end{smallmatrix}\right) \in GL(2,\Rmath)$
such that
$\tau_1 M^{\prime} \tau_1=GMG^{-1}$
with $G$ as in (\ref{G}),
i.e.
$M^{\prime}=(\tau_1 G)M(\tau_1G)^{-1}= G_1 M G_1^{-1}$
with $G_1=\tau_1G=
\left(\begin{smallmatrix}
1&0&0&0\\
0&b^8_8&0&b^8_4\\
0&0&1&0\\
0&b^4_8&0&b^4_4
\end{smallmatrix}\right)$
which is simply the matrix corresponding to
$\left(\begin{smallmatrix}b^8_8&b^8_4\\b^4_8&b^4_4 \end{smallmatrix}\right) \in GL(2,\Rmath)$
in the formula (\ref{G}).
\end{proof}
\begin{corollary}
\label{equivCSu22}
Any $J \in {\frak{X}}_{\mathfrak{u}(2) \oplus \mathfrak{u}(2)}$
is equivalent
under some member of $\left( SO(3)\times \Rmath_+^*\right)^2$
to  $K(M)$ in
(\ref{Ju2xu2notorsion})
with
$M\in GL(4,\Rmath),$
$M^2=-I.$   $K(M),$
$K(M^{\prime})$ are equivalent if and only if
there exists some
$\left(\begin{smallmatrix}b^4_4&b^4_8\\b^8_4&b^8_8 \end{smallmatrix}\right) \in GL(2,\Rmath)$
such that
$M^{\prime}=GMG^{-1}$ with $G$ as in
(\ref{G}).
\end{corollary}
\begin{proof}
Follows readily from
theorem \ref{torsionfreeu22}.
\end{proof}
\section{$U(2)^N.$}
The results of
lemma \ref{automu22},
theorem \ref{torsionfreeu22}
and corollary \ref{equivCSu22}
generalize
in the following way.
\begin{lemma}
For any $N \in \Nmath^*,$
$\Aut{(\mathfrak{u}(2))^N} =
H_N \cup \, \left( \bigcup_{\sigma \in \Sigma} \, \tau_\sigma \, H_N \right)$
(disjoint reunion)
where:
\begin{itemize}
\item
$H_N=\left\{
(U^i_j)_{1\leqslant i,j \leqslant N};
U^i_i=
\left(\begin{smallmatrix}
\Phi_i&0\\0&b^i_i
\end{smallmatrix}
\right),
U^i_j=
\left(\begin{smallmatrix}
0&0\\0&b^j_i
\end{smallmatrix}
\right) ( i \neq j),
\Phi_i \in SO(3), \det{(b^j_i)}\neq 0
\right\};$
\item
$\Sigma$ is the set of circular permutations of $\{1, \cdots, N\}$ having no fixed point,
and $\tau_\sigma=(T^i_j)_{1\leqslant i,j \leqslant N}$
with the $T^i_j$s the $4\times 4$ blocks
$T^i_j= \delta_{i,\sigma(j)}\, I$ ($I$ the $4\times 4$ identity and $\delta_{k,\ell}$ the Kronecker symbol).
\end{itemize}
\end{lemma}
\begin{theorem}
Let  $J \, : \, \mathfrak{u}(2)^N  \rightarrow
\mathfrak{u}(2)^N.$
 $J$ has zero torsion if and only if
there exists
$\Phi \in  \left( SO(3)\times \Rmath_+^*\right)^N \subset H_N$
and $M=(M^i_j)_{1\leqslant i,j \leqslant N} \in \text{gl}(2N,\Rmath),$
$M^i_j=
\begin{pmatrix}
\xi^{4i-1}_{4j-1} &\xi^{4i-1}_{4j}\\
\xi^{4i}_{4j-1} &\xi^{4i}_{4j}
\end{pmatrix},$
such that
$\Phi^{-1}J\Phi = K(M)$
with
\begin{equation}
\label{KMu2N}
K(M) =(K^i_j(M))_{1\leqslant i,j \leqslant N}
\end{equation}
and the $K^i_j(M)$s the following $4\times 4$ blocks:
\begin{multline*}
\label{torsionfreesuN}
K^i_i(M) =    \left(
\begin{array}{cc|c}
0&1&0\\
-1&0&0\\
\hline
0&0&M^i_i
\end{array}\right) \quad (1\leqslant i \leqslant N)\; ,
\\
K^i_j(M) =    \left(
\begin{array}{cc|c}
0&0&0\\
0&0&0\\
\hline
0&0&M^i_j
\end{array}\right)
\quad (1\leqslant i,j \leqslant N, \; i \neq j).
\end{multline*}
\end{theorem}
(Here we used that the analogs of
$17|3,27|3$
$18|3,28|3$
$35|7,36|7$
$45|7,46|7$
at the end of (i) in the proof of Theorem
\ref{torsionfreeu22}
are resp. , with $i <j,$
\begin{eqnarray*}
4i-3,4j-1|4i-1 : &\xi^{4i-3}_{4j-1} -\xi^{4i-1}_{4i-1}\xi^{4i-2}_{4j-1}=0\\
4i-2,4j-1|4i-1 : &\xi^{4i-2}_{4j-1} +\xi^{4i-1}_{4i-1}\xi^{4i-3}_{4j-1}=0\\
4i-3,4j|4i-1 : &\xi^{4i-3}_{4j} -\xi^{4i-1}_{4i-1}\xi^{4i-2}_{4j}=0\\
4i-2,4j|4i-1 : &\xi^{4i-2}_{4j} +\xi^{4i-1}_{4i-1}\xi^{4i-3}_{4j}=0\\
4i-1,4j-3|4j-1 : &\xi^{4j-3}_{4i-1} -\xi^{4j-1}_{4j-1}\xi^{4j-2}_{4i-1}=0\\
4i-1,4j-2|4j-1 : &\xi^{4j-2}_{4i-1} +\xi^{4j-1}_{4j-1}\xi^{4j-3}_{4i-1}=0\\
4i,4j-3|4j-1 : &\xi^{4j-3}_{4i} -\xi^{4j-1}_{4j-1}\xi^{4j-2}_{4i}=0\\
4i,4j-2|4j-1 : &\xi^{4j-2}_{4i} +\xi^{4j-1}_{4j-1}\xi^{4j-3}_{4i}=0
\end{eqnarray*}
and give
$\xi^{4i-3}_{4j-1}=\xi^{4i-2}_{4j-1}=
\xi^{4i-3}_{4j}=\xi^{4i-2}_{4j}=
\xi^{4j-3}_{4i-1}=\xi^{4j-2}_{4i-1}=
\xi^{4j-3}_{4i}=\xi^{4j-2}_{4i}=
0.$ Then all torsion equations vanish.)

\begin{corollary}
Any $J \in {\frak{X}}_{\mathfrak{u}(2)^N}$
is equivalent
under some member of $\left( SO(3)\times \Rmath_+^*\right)^N$
to  $K(M)$ in
(\ref{KMu2N})
with
$M\in GL(2N,\Rmath),$
$M^2=-I.$   $K(M),$
$K(M^{\prime})$ are equivalent if and only if
there exists some
$\left(b^{4i}_{4j}\right)_{1\leqslant i,j \leqslant N}
\in GL(N,\Rmath)$
such that
$M^{\prime}=GMG^{-1}$ with
$G
=(G^i_j(M))_{1\leqslant i,j \leqslant N} ,\;
G^i_i=
\left(\begin{smallmatrix}
1&0\\0&b^{4i}_{4i}
\end{smallmatrix}\right),\;
G^i_j=
\left(\begin{smallmatrix}
0&0\\0&b^{4i}_{4j}
\end{smallmatrix}
\right) ( i \neq j).$
\end{corollary}

\begin{remark}
\rm
The closed set $\mathcal{R}= \{M \in GL(2N,\Rmath); M^2=-I\}$ is comprised
of the conjugates of
$\mathcal{T}=\left(\begin{smallmatrix}
0&-I_N\\
I_N&0
\end{smallmatrix}\right)$ ($I_N$ the $N\times N$ identity)
under the action of $GL(2N,\Rmath).$ Hence it is a $2N^2$-dimensional submanifold
of $\Rmath^{4N^2}$ with a diffeomorphism
$\chi \, : \, GL(2N,\Rmath) \left/ \right. \mathcal{S} \rightarrow \mathcal{R},$
$\mathcal{S} =\left\{ Q=\left(\begin{smallmatrix} R &-S\\ S&R\end{smallmatrix}\right); R,S\in GL(N,\Rmath), \det{Q}\neq 0 \right\}$
 the stabilizer of $\mathcal{T},$ and $\chi$
defined by $\chi\,  \left[ P\right]\, = P\mathcal{T}P^{-1} $ for
$\left[ P\right]$ the class mod
$\mathcal{S}$ of  $P \in GL(2N,\Rmath).$
For $N=2,$
$\chi \left[\left( \begin{smallmatrix} -\eta&0\\ \xi&1 \end{smallmatrix} \right) \right]
= \left( \begin{smallmatrix} \xi&\eta\\-\frac{1+\xi^2}{\eta}&- \xi \end{smallmatrix} \right),$ $(\xi,\eta)\in \Rmath \times \Rmath^*.$
For general $N$ and for $G
\in GL(2N,\Rmath), M=\chi \left[P\right],
M^{\prime} =
\chi\, \left[P^{\prime}\right] \in \mathcal{R},$
$GMG^{-1}=
\chi\, \left[GP\right]\,$ and the condition
$M^{\prime}=GMG^{-1}$
 reads
$\left[P^{\prime}\right]=
 \left[GP\right].$

\end{remark}

\section{Local chart and a representation for $(U(2), J(\xi)).$}
\subsection{Local chart.}
For any fixed $\xi \in \Rmath,$
denote simply $J$ the complex structure
$J(\xi)$ on
$\mathfrak{u}(2)$
and by $G$ the group $U(2)$ endowed with the left invariant structure
of complex manifold  defined by $J$.
For any open subset $V  \subset U(2),$ the space $H_{\Cmath}(V)$
of complex valued holomorphic functions on $V$ (considered here as a subset of $G$)
is comprised of all complex smooth functions $f$ on $V$
 which are
annihilated by
all
$\tilde{X}_j^{-} = X_j +i J X_j$ ($ 1 \leqslant j \leqslant 4$)
with $(X_j)_{1 \leqslant j \leqslant 4}$
(resp. $(JX_j)$)
the left invariant vector fields associated
to the basis $(J_j)_{1 \leqslant j \leqslant 4}$ of
$\mathfrak{u}(2)$
(resp. to $(J\,  J_j)$).
One has
$\tilde{X}_1^{-} = X_1 -i X_2,$
$\tilde{X}_2^{-} = i\tilde{X}_1^{-},$
$\tilde{X}_4^{-} = iX_3 +(1-i\xi) X_4,$
$\tilde{X}_3^{-} = - i(1+i\xi)\tilde{X}_4^{-},$
hence
\begin{equation}
\label{H(V)}
H_{\Cmath}(V) =\{f \in C^{\infty}(V) \; ; \;
\tilde{X}_1^{-} \, f =
\tilde{X}_4^{-} \, f = 0
\}.
\end{equation}
As is known, the map $\Smath^1 \times SU(2) \rightarrow U(2)$ defined by
$$(\zeta, A) \mapsto
\begin{pmatrix}
\zeta&0\\0&1
\end{pmatrix}
A$$
is a diffeomorphism of manifolds (not of groups).
Introducing Euler angles as coordinates in the open subset
$\Omega =SU(2) \setminus \left(
e^{\Rmath J_1} \cup
e^{\pi J_3} e^{\Rmath J_1} \right)$
of $SU(2),$  one gets the coordinates
$(s,\theta,\varphi,\psi)$ in the open subset
\begin{equation}
\label{V}
V= \left( \Smath^1 \setminus \{-1\} \right) \times \Omega
\end{equation}
such that $u$ defined by
\begin{equation}
\label{xangles euler}
u(s,\theta,\varphi,\psi) =
\begin{pmatrix}
e^{is}&0\\0&1
\end{pmatrix}
e^{\varphi J_3}
e^{\theta J_1}
e^{\psi J_3} =
\begin{pmatrix}
e^{is}e^{i\frac{\varphi+\psi}{2}}\cos{\frac{\theta}{2}}&
ie^{is}e^{i\frac{\varphi-\psi}{2}}\sin{\frac{\theta}{2}}\\
ie^{-i\frac{\varphi-\psi}{2}}\sin{\frac{\theta}{2}}&
e^{-i\frac{\varphi+\psi}{2}}\cos{\frac{\theta}{2}}
\end{pmatrix}
\end{equation}
is a diffeomorphism of $
]-\pi,\pi[ \times
]0,\pi[ \times
]0,2\pi[ \times
]-2\pi,2\pi[ $ on $V.$
Then one gets on $V$ (see e.g. \cite{Vilenkin}, p.141):
\begin{eqnarray*}
X_1 &=& \cos{\psi}\, \frac{\partial}{\partial \theta}
+\frac{\sin{\psi}}{\sin{\theta}}\, \frac{\partial}{\partial \varphi}
-\cot{\theta}\sin{\psi}\, \frac{\partial}{\partial \psi}
\\
X_2 &=& -\sin{\psi}\, \frac{\partial}{\partial \theta}
+\frac{\cos{\psi}}{\sin{\theta}}\, \frac{\partial}{\partial \varphi}
-\cot{\theta}\cos{\psi}\, \frac{\partial}{\partial \psi}
\\
X_3 &=&  \frac{\partial}{\partial \psi} \\
X_4 &=&  \frac{\partial}{\partial s} - \frac{\partial}{\partial \varphi}
\end{eqnarray*}
Hence
$f \in C^{\infty}(V)$ is in $ H_{\Cmath}(V)$ if and only if it
satisfies the 2 equations
\begin{eqnarray}
i\, {\sin{\theta}}\, \frac{\partial f}{\partial \theta}
+  \frac{\partial f }{\partial \varphi}
-\cos{\theta}\, \frac{\partial f}{\partial \psi}  &=&0\\
i\, \frac{\partial f}{\partial \psi}
+ (1-i\xi)\left( \frac{\partial f}{\partial s}
                  - \frac{\partial f}{\partial \varphi}\right) &=&0
\end{eqnarray}
The 2 functions
\begin{eqnarray}
\label{w1}
w^1 &=&e^{i(s+\varphi)} \cot{\frac{\theta}{2}}\\
\label{w2}
w^2&=& e^{(1+i\xi)\frac{s}{2(1+\xi^2)}} e^{i\frac{\psi}{2}}\sqrt{\sin{\theta}}
\end{eqnarray}
are holomorphic on $V.$
Let $F : V \rightarrow \Cmath^2$ defined
by $F=(w^1,w^2).$ It is easily seen that $F$ is injective, with jacobian
$-\frac{1}{4(1+\xi^2)}\, e^{\frac{s}{1+\xi^2}}\, (\cot{\frac{\theta}{2}})^2 \neq 0,$
hence
$F$ is a biholomorphic bijection of $V$ onto an open subset $F(V)$ of $\Cmath^2,$
i.e. $(V,F)$ is a  chart of $G.$
$F(V)$ is the set of those
$(w^1,w^2)\in \Cmath^2$ satisfying
the following conditions, where $r_1=|w^1|,r_2=|w^2|$ and
$\omega(r_1,r_2)=\log{r_2}-\frac{1}{2}\,\log{\frac{2r_1}{1+r_1^2}}:$
$r_1r_2 \neq 0,$ ,
${\sqrt{\frac{2r_1}{1+r_1^2}}}  \, e^{-\frac{\pi}{2(1+\xi^2)}}
 < r_2
< {\sqrt{\frac{2r_1}{1+r_1^2}}}  \, e^{\frac{\pi}{2(1+\xi^2)}}
$ ,
$\arg{w_1} \not \equiv 2
(1+\xi^2)
\omega(r_1,r_2)
\, \text{ mod } 2\pi$ ,
$\arg{w_2} \not \equiv \xi \omega(r_1,r_2) +\pi \, \text{ mod } 2\pi$.
For example, if $\xi=0,$
$$
V=
\bigcup_{r_1>0}\quad
\bigcup_{e^{-\frac{\pi}{2}}y(r_1)< r_2  < e^{\frac{\pi}{2}}y(r_1)}
\, \left(
\left(
\mathcal{C}^{(1)}_{r_1} \setminus \left\{ \text{arg } \equiv 2\omega(r_1,r_2) \right\}
\right)
\times
\left(
\mathcal{C}^{(2)}_{r_2} \setminus \left\{ \text{arg } \equiv \pi \right\}
\right)
\right)
$$
where $\mathcal{C}^j_{r_j}$ ($j=1,2$) is the circle with radius $r_j$ in the $w^j$-plane
and $y(x)= \sqrt{\frac{2x}{1+x^2}}$ $(x>0).$

\subsection{A representation on a space of holomorphic functions.}
As $U(2)$ is compact, there are no nonconstant holomorphic functions on the whole
of $U(2).$ Instead, we consider  the space $H_{\Cmath}(V)$
of holomorphic functions on the open subset $V$
(\ref{V}), and we  compute
(as kind of substitute for the regular representation)
the representation $\lambda$
of the Lie algebra
$\mathfrak{u}(2)$
we get by Lie derivatives on
 $H_{\Cmath}(V).$
 First, note that for any $x =  \left(\begin{smallmatrix} a&b\\c&d\end{smallmatrix} \right) \in V$
 as in (\ref{xangles euler}), the  complex coordinates
 $w^1,w^2$ of $x$
(\ref{w1},\ref{w2})  satisfy:
 \begin{eqnarray}
 w^1&=&-i e^{-is} \frac{a}{\bar{b}}\\
 \left(w^2\right)^2&=&2i a \bar{b}\, e^{s \frac{1+i\xi}{1+\xi^2}}.
 \end{eqnarray}
 Then one gets for the complex coordinates
 $w^1_{e^{-tJ_1}x}, w^2_{e^{-tJ_1}x} $ of
 $e^{-tJ_1}x$ ($x \in V,$ $t\in \Rmath$ sufficiently small):
 \begin{eqnarray*}
 w^1_{e^{-tJ_1}x}&=&\frac{1+w^1\, \cot{\frac{t}{2}} }{\cot{\frac{t}{2}}-w^1}\\
 (w^2_{e^{-tJ_1}x})^2&=&\left(w^2\right)^2 \left( \cos{t} +
 \frac{\sin{t}}{2} \, \frac{1- \left(w^1\right)^2}{w^1}\right).
 \end{eqnarray*}
 Whence for any  $f  \in H_{\Cmath}(V)$, denoting $J_1f$ instead of
 $\lambda(J_1)f:$
 \begin{equation*}
 \left(J_1f\right) (w^1,w^2)  = \left[\frac{d}{dt} f(w^1_{e^{-tJ_1}x}, w^2_{e^{-tJ_1}x} ) \right]_{t=0}
 =
 \frac{1+(w^1)^2}{2}  \,
 \frac{\partial f}{\partial w^1}
 +
 \frac{w^2\left(1-(w^1)^2\right)}{4w^1}
 \, \frac{\partial f}{\partial w^2}.
 \end{equation*}
 In the same way,
 \begin{eqnarray*}
 w^1_{e^{-tJ_2}x}&=&-i \, \frac{i \sin{\frac{t}{2}} + w^1 \, \cos{\frac{t}{2}} }
 {-i\cos{\frac{t}{2}} + w^1\, \sin{\frac{t}{2}} }
 \\
 (w^2_{e^{-tJ_2}x})^2&=&\left(w^2\right)^2 \left( \cos{t} +
 i \, \frac{\sin{t}}{2} \, \frac{1+ \left(w^1\right)^2}{w^1}\right)
 \\
 \left(J_2f\right) (w^1,w^2)
 &=&
 \frac{i(1-(w^1)^2)}{2}  \,
 \frac{\partial f}{\partial w^1}
 +
 \frac{iw^2\left(1+(w^1)^2\right)}{4w^1}
 \, \frac{\partial f}{\partial w^2}
 \\
 w^1_{e^{-tJ_3}x}&=& e^{-it} w^1 \\
 (w^2_{e^{-tJ_3}x})^2&=&\left(w^2\right)^2 \\
 \left(J_3f\right) (w^1,w^2)
 &=&      -i w^1
 \frac{\partial f}{\partial w^1}.
 \end{eqnarray*}
 Finally,
 \begin{eqnarray*}
 w^1_{e^{-tJ_4}x}&=&w^1\\
 (w^2_{e^{-tJ_4}x})^2&=&\left(w^2\right)^2 e^{-t\frac{1+i\xi}{1+\xi^2}} \\
 \left(J_4f\right) (w^1,w^2)
 &=&      -  \frac{1+i\xi}{1+\xi^2}  \,
 w^2 \,
 \frac{\partial f}{\partial w^2}.
 \end{eqnarray*}
In the complexification $\mathfrak{sl}(2) \oplus \Cmath J_4   $ of $\mathfrak{u}(2),$
introduce as usual
 \begin{equation*}
 H_{\pm} =iJ_1 \mp J_2, \; H_3=iJ_3
 \end{equation*}
 so that
 \begin{equation*}
 [H_3,H_{\pm}] =\pm H_{\pm} , \; [H_+,H_-] = 2H_3.
 \end{equation*}
 Then, extending the representation $\lambda$ to
$\mathfrak{sl}(2) \oplus \Cmath J_4,$
one has, with
$H_4 =      -  (1-i\xi) J_4$
:
 \begin{eqnarray}
 \left(H_+f\right) (w^1,w^2)
 &=& i \left(
 (w^1)^2
 \frac{\partial f}{\partial w^1}
 -
 \frac{1}{2}
 \, w^1w^2 \,\frac{\partial f}{\partial w^2}
 \right)
 \\
 \left(H_- f\right) (w^1,w^2)
 &=& i \left(
 \frac{\partial f}{\partial w^1}
 +
 \frac{w^2}{2w^1}
 \,\frac{\partial f}{\partial w^2}
 \right)
 \\
 \left(H_3f\right) (w^1,w^2)
 &=&       w^1
 \frac{\partial f}{\partial w^1}.
 \\
 \left(H_4f\right) (w^1,w^2)
 &=&
 w^2 \,
 \frac{\partial f}{\partial w^2}.
 \end{eqnarray}

\subsection{A subrepresentation.}
We restrict $\lambda$ to $H(\Cmath^* \times \Cmath^*)$
($\Cmath^*= \Cmath \setminus \{0\}$),
and denote $\varphi_{p,q}$ the function
$\varphi_{p,q}(w^1,w^2)= (w^1)^p(w^2)^q$ for $p,q \in \Zmath.$
The system $(\varphi_{p,q})_{p,q \in \Zmath}$ is total in
$H(\Cmath^* \times \Cmath^*),$ and one has:
 \begin{eqnarray}
 \label{repeq1}
 H_+ \,
 \varphi_{p,q}
 &=& i \left(
 p- \frac{q}{2}\right)
 \varphi_{p+1,q}
 \\
 \label{repeq2}
 H_-   \,
 \varphi_{p,q}
 &=& i \left(
 p+ \frac{q}{2}\right)
 \varphi_{p-1,q}
 \\
 H_3\,
 \varphi_{p,q}
 &=& p\,
 \varphi_{p,q}
 \\
 H_4    \,
 \varphi_{p,q}
 &=& q    \,
 \varphi_{p,q}
 \end{eqnarray}
 For any $q \in \Zmath,$
 the subspace $\mathcal{H}_q$ of functions of the form
 $(w^2)^q g(w^1),$ $g \in H(\Cmath^*),$ is a closed invariant
 subspace
 of
$H(\Cmath^* \times \Cmath^*),$ and
$H(\Cmath^* \times \Cmath^*)$ is the closure of
$\bigoplus_{q\in \Zmath}  \mathcal{H}_q.$

\subsection{A lemma.}
\begin{lemma}
\label{topirredu}
Let
$\mathcal{E} = H(\Cmath^*)$
the Fr\'{e}chet space
of holomorphic functions of the complex variable $z$
on $\Cmath^*.$
Let $\mathcal{F}$ be any closed vector subspace
of
$\mathcal{E}$ that is invariant by the  operator
$z \frac{d }{d z}.$
Let $ f \in \mathcal{F}$ and
$f(z) = \sum_{p=-\infty}^{+\infty} \; c_p z^p$ its Laurent expansion in $\Cmath^*.$
If for some $p \in \Zmath, c_p \neq 0,$ then the function  $z \mapsto z^p$ belongs to
$\mathcal{F}.$
\end{lemma}
\begin{proof}
We show first that
the function $z \mapsto f(e^{i\theta}z)$ belongs to
$\mathcal{F}$ $\forall \theta \in \Rmath$
 $\forall f \in \mathcal{F}.$
Let $ f \in \mathcal{F}$ and
$f(z) = \sum_{p=-\infty}^{+\infty} \; c_p z^p$ its Laurent expansion in $\Cmath^*.$
Since it is uniformly and absolutely convergent on compact subsets of
$\Cmath^*,$ and since the operator $H= z \frac{d}{dz}$ is continuous on
$\mathcal{E},$
$$\frac{(i\theta)^k}{k!}\,  (H^kf)(z)=
\sum_{p=-\infty}^{\infty} \; c_p \frac{(i\theta p)^k}{k!} \, z^p, \quad \forall k \in \Nmath,
\, \forall \theta \in \Rmath,
\, \forall z \in \Cmath^*.
$$
On the other hand, for any fixed $\theta\in \Rmath,$ the double series
$$\sum_{k=0}^{+\infty} \; \sum_{p=-\infty}^{+\infty} \; c_p \frac{(i\theta p)^k}{k!}z^p$$
is absolutely and uniformly summable in the annulus $A(r,R)$ for any $0<r< R<+\infty$
since
$$\sum_{k=0}^{+\infty} \; \sum_{p=-\infty}^{+\infty} \; |c_p| \frac{(|\theta| |p|)^k}{k!}
|z|^p \leqslant
\sum_{p<0} \; |c_p| (e^{-|\theta|}r)^p
+
\sum_{p\geqslant 0} \; |c_p| (e^{|\theta|}R)^p < +\infty.$$
From the associativity theorem for summable families,
$$f(e^{i\theta}z) =
\sum_{k=0}^{+\infty} \;
\frac{(i\theta)^k}{k!}\,  (H^kf)(z)$$
with the series
uniformly and absolutely convergent on compact subsets of
$\Cmath^*.$
The conclusion follows, since
$H^kf \in \mathcal{F}$ $\forall k.$
Now we use the same trick as in
(\cite{helgason}, p. 14). For any $z \in \Cmath^*,$ denote $f_z$ the periodic function on
$\Rmath:$ $\theta \mapsto f(e^{i\theta} z).$ Its Fourier expansion is
$f(e^{i\theta} z) = \sum_{p =-\infty}^{+\infty} \tilde{c}_p(z) e^{i p \theta}$
where
$$\tilde{c}_p(z) = \frac{1}{2\pi} \int_{0}^{2\pi}
f(e^{i\theta} z) e^{-i p \theta} d\theta.$$
The function $z \mapsto \tilde{c}_p(z)$ belongs to
$\mathcal{F}$ as the right-hand side is a limit in
$\mathcal{E}$ of linear combinations of functions $z \mapsto f(e^{i\theta} z)$.
But with the Laurent expansion of $f$ one gets
$f(e^{i\theta} z) = \sum_{p= -\infty}^{+\infty} c_p z^p e^{i p \theta}.$
For any $z,$
that series is a trigonometric series that converges uniformly on $\Rmath$
hence it coincides with the Fourier series of $f_z$
and
$\tilde{c}_p(z)= c_p z^p $ $\forall p \in \Zmath.$
Hence if for some $p \in \Zmath, c_p \neq 0,$
then the function  $z \mapsto z^p$ belongs to
$\mathcal{F}.$
\end{proof}

\subsection{A closer look to the subrepresentation.}
Introduce the Casimir $C=H_+H_- +\left(H_3\right)^2 -H_3.$
On
 $\mathcal{H}_q,$  $C= u(u+1)$ with $u=\frac{q}{2}.$
 Now, we distinguish cases.
 We use both the notations
 $\uparrow_{u}^{q},
 D^{q}(2k),$ etc.
 of \cite{miller} (Th. 2.3)
 for representations of $\frak{u}(2)$
 and usual notations of e.g.
  \cite{miller2} (7.3)
 $\uparrow_{u},
 D^{(k)},$  etc.
 for representations of
 $\frak{sl}(2).$  One has
 $\uparrow_{u}^{q} =  \uparrow_{u} \otimes \, q$  , \,
 $D^q(2k) =  D^{(k)} \otimes \, q$ etc.
 \\ \textbf{Case 1:}
 $q=-2k,$ $k\in \Nmath \setminus \{0\}.$
 Then
 from
 (\ref{repeq1}), (\ref{repeq2}),
 the closed subspace
 $\mathcal{H}_{q}^{\uparrow}$
 (resp. $\mathcal{H}_{q}^{\downarrow}$)
 generated by
 $\{ \varphi_{k+n,q}, n\in \Nmath\},$
 (resp. $\{ \varphi_{-k-n,q}, n\in \Nmath\}$),
 which is comprised of the functions $(w^2)^{-2k} (w^1)^{k} g(w^1)$
 (resp. $(w^2)^{-2k} (w^1)^{-k} g(\frac{1}{w^1})$),
 $g \in H(\Cmath),$ is invariant and topologically irreducible from
 lemma \ref{topirredu}.
 $\mathcal{H}_{q}^{\uparrow} = \uparrow_{-k}^{q} =  \uparrow_{-k} \otimes \, q$
 , \,
 $\mathcal{H}_{q}^{\downarrow} = \downarrow_{-k}^{q} =  \downarrow_{-k} \otimes \, q$
 .
 $\mathcal{H}_{q}$ is indecomposable and
 $\mathcal{H}_{q}\left/ \left(\mathcal{H}_{q}^{\uparrow} \oplus  \mathcal{H}_{q}^{\downarrow}\right)
 \right.$
 is top. irreducible and equal to $D^{q}(2(k-1)) = D^{(k-1)} \otimes q
 ,$ i.e.
 $\mathcal{H}_{q}$ is a nontrivial extension of
 $D^{q}(2(k-1))$ by
 $\uparrow_{-k}^{q} \oplus  \downarrow_{-k}^{q}.$
 \\ \textbf{Case 2:}
 $q=2k,$ $k\in \Nmath.$
 Then the closed subspace
 $\mathcal{H}_{q}^{D}$
 generated by
 $\{ \varphi_{-k+n,q}, n\in \Nmath, 0\leqslant n \leqslant 2k\},$
 which is comprised of the functions $(w^2)^{2k} (w^1)^{-k} P(w^1),$
 $P \in \Cmath[w^1], \deg{P}\leqslant 2k,$ is invariant and topologically irreducible from
 lemma \ref{topirredu}, and
 $\mathcal{H}_{q}^{D} = D^{q}(2k)=   D^{(k)} \otimes q.$
 There are exactly 2 closed invariant (nontrivial) subspaces containing
 $\mathcal{H}_{q}^{D}.$ Each one is indecomposable, with top. irreducible quotient
 by  $\mathcal{H}_{q}^{D}$ equal respectively to
 $\uparrow_{-k-1}^{q}$ or $\downarrow_{-k-1}^{q}.$
 \\ \textbf{Case 3:} $q\not \in 2\Zmath.$
 In that case
 $\mathcal{H}_{q} = D^{q}(u,0)$ is top. irreducible.

 We see that $\lambda$ is quite different from the regular representation,
 since the differentials of the representation in the  unitary dual of $U(2)$ are
 $D^{(\ell)} \otimes \, m$,  $2 \ell \in \Nmath,  m \in \Zmath,$ with $2\ell +m $ even
 (\cite{brocker-dieck}, p. 87).
\section{Chart for $(SU(2) \times SU(2), J(\xi,\eta)).$}
In this last section, we compute an holomorphic chart for
$J(\xi,\eta),$
$(\xi, \eta) \in \Rmath\times \Rmath^*,$
in the open subset $W = \Omega\times \Omega$ of
$SU(2)\times SU(2)$  with Euler angles coordinates
$(\theta_1,\phi_1,\psi_1,\theta_2,\phi_2,\psi_2).$
The space  $H_{\Cmath}(W)$  of complex valued holomorphic functions on $W$
is comprised of all complex smooth functions $f$ on $W$
 which are annihilated by
all
\begin{equation}
{ \left.\tilde{X}_j^{(k)}\right. }^{-} = X_j^{(k)} +i J X_j^{(k)},
 \quad 1 \leqslant j \leqslant 3,
 \quad 1 \leqslant k \leqslant 2,
\end{equation}
$(X_j^{(k)})$ the left invariant vector fields associated
to the basis
$( J_1^{(1)},
 J_2^{(1)},
 J_3^{(1)},
 J_1^{(2)},
 J_2^{(2)},
 J_3^{(2)})$
of $\mathfrak{su}(2) \oplus \mathfrak{su}(2).$
One has
${ \left.\tilde{X}_1^{(k)}\right. }^{-} = X_1^{(k)} -i X_2^{(k)},$
${ \left.\tilde{X}_2^{(k)}\right. }^{-} = i{ \left.\tilde{X}_1^{(k)}\right. }^{-} $
for $(k=1,2),$ and
${ \left.\tilde{X}_3^{(2)}\right. }^{-} = i \eta X_3^{(1)} +(1-i\xi) X_3^{(2)},$
${ \left.\tilde{X}_3^{(1)}\right. }^{-} =
-i\frac{1+i\xi}{\eta}{ \left.\tilde{X}_3^{(2)}\right. }^{-}.$
For $k=1,2,$
\begin{eqnarray*}
X_1^{(k)} &=& \cos{\psi_k}\, \frac{\partial}{\partial \theta_k}
+\frac{\sin{\psi_k}}{\sin{\theta_k}}\, \frac{\partial}{\partial \varphi_k}
-\cot{\theta_k}\sin{\psi_k}\, \frac{\partial}{\partial \psi_k}
\\
X_2^{(k)} &=& -\sin{\psi_k}\, \frac{\partial}{\partial \theta_k}
+\frac{\cos{\psi_k}}{\sin{\theta_k}}\, \frac{\partial}{\partial \varphi_k}
-\cot{\theta_k}\cos{\psi_k}\, \frac{\partial}{\partial \psi_k}
\\
X_3^{(k)} &=&  \frac{\partial}{\partial \psi_k}
\, .
\end{eqnarray*}
Hence
$f \in C^{\infty}(W)$ is in $ H_{\Cmath}(W)$ if and only if it
satisfies the 3 equations
\begin{eqnarray}
i\, {\sin{\theta_1}}\, \frac{\partial f}{\partial \theta_1}
+  \frac{\partial f }{\partial \varphi_1}
-\cos{\theta_1}\, \frac{\partial f}{\partial \psi_1}  &=&0\\
i\, {\sin{\theta_2}}\, \frac{\partial f}{\partial \theta_2}
+  \frac{\partial f }{\partial \varphi_2}
-\cos{\theta_2}\, \frac{\partial f}{\partial \psi_2}  &=&0\\
i\eta\, \frac{\partial f}{\partial \psi_1}
+ (1-i\xi) \frac{\partial f}{\partial \psi_2}
&=&0 \, .
\end{eqnarray}
The 3 functions
\begin{eqnarray}
\label{z1}
z^1 &=&e^{i\varphi_1} \cot{\frac{\theta_1}{2}}\\
\label{z2}
z^2 &=&e^{i\varphi_2} \cot{\frac{\theta_2}{2}}\\
\label{z3}
z^3&=&
e^{i\frac{\psi_1}{2}}
e^{\frac{\eta(1+i\xi)}{1+\xi^2}\,\frac{ \psi_2}{2} }
\sqrt{\sin{\theta_1}}
\sqrt{\sin{\theta_2}}
\end{eqnarray}
are holomorphic on $W.$
Let $Z : W \rightarrow \Cmath^3$ defined
by $Z=(z^1,z^2,z^3).$  $Z$ is injective, with jacobian
$
-\frac{\eta}{4(1+\xi^2)} \, e^{\frac{\eta}{1+\xi^2}\, \psi_2}
(\cot{\frac{\theta_1}{2}})^2
(\cot{\frac{\theta_2}{2}})^2
\neq 0,$
hence
$Z$ is a biholomorphic bijection of $W$ onto an open subset  of $\Cmath^3,$
i.e. $(W,Z)$ is a  local chart for  $SU(2)\times SU(2)$
equipped with the complex structure $J(\xi,\eta).$



\begin{thebibliography}{99}

\bibitem{brocker-dieck}
 T. Br\"{o}cker, T. Dieck
\textit{Representations of Compact Lie Groups},
Graduate texts in Math. \#98,
Springer, New York, 1985.


\bibitem{charbonnel}
J.-Y. Charbonnel,  H. O. Khalgui,
Classification des structures $CR$ invariantes pour les groupes de Lie compacts,
\textit{J. of Lie Theory}, \textbf{14}, 2004,  165-198.

\bibitem{companionarchive}
\texttt{http://www.u-bourgogne.fr/monge/l.magnin/CSu2/CSu2index.html}
or \texttt{http://math.u-bourgogne.fr/IMB/magnin/public\_html/CSu2/CSu2index.html}



\bibitem{daurtseva}
N. A. Daurtseva,
Invariant complex structures on $\Smath^3 \times \Smath^3$,
Electronic journal
\textit{Investigated in Russia},
 2004,   888-893.\\
English version
\texttt{http://zhurnal.ape.relarn.ru/articles/2004/081e.pdf}
\\
Russian version
\texttt{http://zhurnal.ape.relarn.ru/articles/2004/081.pdf}


\bibitem{helgason}
 S. Helgason,
\textit{Groups and geometric analysis (integral geometry, invariant differential operators
and spherical functions)},
Academic Press, Orlando, 1984.

\bibitem{artmagnin1}
L. Magnin, Complex structures on indecomposable
     6-dimensional nilpotent real Lie algebras,
\textit{Int. J. Alg. Comput.},
\textbf{17}, \#1, 2007,  77-113.

\bibitem{artmagnin2}
L. Magnin, Left invariant complex structures on
real 6-dimensional simply connected indecomposable nilpotent Lie groups,
\textit{Int. J. Alg. Comput.},
\textbf{17}, \#1, 2007,  115-139.



\bibitem{miller}
 W.  S. Miller Jr.,
\textit{Lie theory and special functions},
Academic Press, New York, 1968.

\bibitem{miller2}
 W.  Miller Jr.,
\textit{Symmetry groups and their applications},
Academic Press, New York, 1972.



\bibitem{NN}
A. Newlander,  L. Nirenberg,
Complex analytic coordinates in almost complex manifolds,
\textit{Ann. Math.}, \textbf{65}, 1957,  391-404.

\bibitem{sasaki}
T. Sasaki, Classification of left invariant complex structures on $GL(2,\Rmath)$ and $U(2),$
\textit{Kumamoto  J.  Sci. (Math)},
\textbf{14}, 1981,  115-123.

\bibitem{Vilenkin}
 N. Ja. Vilenkin,
\textit{Fonctions sp\'{e}ciales et th\'{e}orie de la repr\'{e}sentation des groupes},
Dunod, Paris, 1969.

\end{thebibliography}
\end{document}